\documentclass[11pt,twoside]{article}
\title{} \author{} \date{}
\usepackage{latexsym,amssymb,times}
\input amssym.def
\newtheorem{te}{Theorem}[section]
\newtheorem{prop}[te]{Proposition}
\newtheorem{cor}[te]{Corollary}

\newtheorem{fac}[te]{Fact}

\newtheorem{ex}[te]{Example}

\def\dok{\noindent{\bf Proof. }}
\def\kdok{\hfill $\Box$ \par \vspace*{2mm} }
\def\a{\alpha}
\def\b{\beta}
\def\g{\gamma}

\def\f{\varphi}
\def\o{\omega}
\def\k{\kappa}

\def\s{\sigma}
\def\p{\psi}
\def\l{\lambda}
\def\t{{\mathfrak t}}
\def\h{{\mathfrak h}}
\def\fb{{\mathfrak b}}

\def\fh{{\mathfrak h}}
\def\ft{{\mathfrak t}}
\def\fc{{\mathfrak c}}
\def\fd{{\mathfrak d}}
\def\fu{{\mathfrak u}}
\def\fp{{\mathfrak p}}

\def\B{{\mathbb B}}
\def\C{{\mathbb C}}

\def\P{{\mathbb P}}
\def\Q{{\mathbb Q}}
\def\R{{\mathbb R}}
\def\B{{\mathbb B}}
\def\N{{\mathbb N}}

\def\BT{{\mathbb T}}

\def\CS{{\mathcal S}}
\def\I{{\mathcal I}}

\def\CP{{\mathcal P}}
\def\CB{{\mathcal B}}
\def\CD{{\mathcal D}}
\def\CM{{\mathcal M}}

\def\CU{{\mathcal U}}

\def\CO{{\mathcal O}}
\def\CA{{\mathcal A}}
\def\CX{{\mathcal X}}
\def\la{\langle}
\def\ra{\rangle}

\def\da{\!\downarrow}

\def\Lim{\mathop{\rm Lim}\nolimits}
\def\Fin{\mathop{\rm Fin}\nolimits}
\def\sm{\mathop{\rm sm}\nolimits}
\def\sq{\mathop{\rm sq}\nolimits}
\def\asq{\mathop{\rm asq}\nolimits}
\def\rp{\mathop{\rm rp}\nolimits}
\def\supp{\mathop{\rm supp}\nolimits}

\def\dom{\mathop{\rm dom}\nolimits}
\def\cf{\mathop{\rm cf}\nolimits}

\def\cc{\mathop{\rm cc}\nolimits}
\def\Col{\mathop{\rm Col}\nolimits}

\def\Emb{\mathop{\rm Emb}\nolimits}

\def\ro{\mathop{\rm ro}\nolimits}

\def\add{\mathop{\mathrm{add}}\nolimits}

\def\Card{\mathop{\mathrm{Card}}\nolimits}

\def\Fn{\mathop{\mathrm{Fn}}\nolimits}
\def\ran{\mathop{\mathrm{ran}}\nolimits}

\def\Borel{\mathop{\mathrm{Borel}}\nolimits}
\def\RO{\mathop{\mathrm{RO}}\nolimits}

\def\Clop{\mathop{\mathrm{Clop}}\nolimits}

\begin{document}
\thispagestyle{plain}
\begin{center}
           {\large \bf {\uppercase{Reduced products of collapsing algebras}}}
\end{center}
\begin{center}
{\bf Milo\v s S.\ Kurili\'c}\footnote{Department of Mathematics and Informatics, Faculty of Sciences, University of Novi Sad,
                                             Trg Dositeja Obradovi\'ca 4, 21000 Novi Sad, Serbia.
                                             e-mail: milos@dmi.uns.ac.rs}
\end{center}
\begin{abstract}
\noindent
$\mathop{\rm rp}\nolimits ({\mathbb B})$ denotes the reduced power ${\mathbb B}^\omega /\Phi$ of a Boolean algebra ${\mathbb B}$, where $\Phi$ is the Fr\'{e}chet filter $\Phi$ on $\omega$.
We investigate iterated reduced powers ($\mathop{\rm rp}\nolimits ^0 ({\mathbb B})={\mathbb B}$
and $\mathop{\rm rp}\nolimits ^{n+1} ({\mathbb B} )=\mathop{\rm rp}\nolimits (\mathop{\rm rp}\nolimits ^n ({\mathbb B}))$) of collapsing algebras
and our main intention is to classify the algebras $\mathop{\rm rp}\nolimits ^n (\mathop{\rm Col}\nolimits (\lambda ,\kappa))$, $n\in {\mathbb N}$, up to isomorphism of their Boolean completions.
In particular, assuming that SCH and ${\mathfrak h} =\omega _1$ hold,
we show that for any cardinals $\lambda\geq \omega$ and $\kappa \geq 2$ such that $\kappa\lambda >\o$ and $\mathop{\rm cf}\nolimits (\lambda )\leq {\mathfrak c}$
we have $\ro (\mathop{\rm rp}\nolimits ^n(\mathop{\rm Col}\nolimits (\lambda ,\kappa)))\cong \mathop{\rm Col}\nolimits (\omega _1, (\kappa ^{<\lambda })^\omega )$, for each $n\in {\mathbb N}$;
more precisely,
$$
\mathop{\rm ro}\nolimits (\mathop{\rm rp}\nolimits ^n(\mathop{\rm Col}\nolimits (\lambda,\kappa)))\cong   \left\{
                          \begin{array}{ll}
                             \mathop{\rm Col}\nolimits (\omega _1,{\mathfrak c})                    ,     & \mbox{ if } \kappa ^{<\lambda} \leq {\mathfrak c} ;\\[1mm]
                             \mathop{\rm Col}\nolimits (\omega _1,\kappa ^{<\lambda}),    & \mbox{ if } \kappa ^{<\lambda} > {\mathfrak c} \land \mathop{\rm cf}\nolimits(\kappa^{<\lambda})>\omega;\\[1mm]
                             \mathop{\rm Col}\nolimits (\omega _1,(\kappa ^{<\lambda}) ^+ ),     & \mbox{ if } \kappa ^{<\lambda} > {\mathfrak c} \land \mathop{\rm cf}\nolimits(\kappa^{<\lambda})=\omega.
                          \end{array}
                                  \right.
$$
If ${\mathfrak b} ={\mathfrak d}$ and $0^\sharp$ does not exist, then the same holds whenever $\mathop{\rm cf}\nolimits (\lambda )= \omega$.
\end{abstract}
\section{Introduction}\label{S1}
In this article we investigate countable powers of Boolean algebras reduced by the Fr\'{e}chet filter $\Phi$ on $\o$.
For a Boolean algebra $\B$, instead of $\B^\o/\Phi$ we write $\rp (\B)$,
and consider iterated reduced powers: $\rp ^0 (\B ):=\B $ and $\rp ^{n+1} (\B ):=\rp (\rp ^n (\B ))$.

The initial motivation for the analysis of such iterations
origins from the investigation of monoids of self embeddings and partial orders of copies of relational structures \cite{Ktow,Kscatt,Kord,Kordu,KurTod,KTR}.
When well orders are in question and $\a \geq \o$ is an ordinal,
the monoid $\Emb (\a)$ of self-embeddings of $\a$ and right Green's preorder on $\Emb (\a )$
are completely reconstructed by the poset $\la \P (\a), \subset\ra$,
where $\P (\a ):=\{ f[\a]:f\in\Emb (\a )\}$ is the set of copies of $\a$ inside $\a$.
So,  in order to obtain a classification of countable ordinals
related to the similarities of their self-embedding monoids and posets of copies,
it was shown in \cite{Kord} that for the separative quotient of the poset $\P (\a)$ we have
\begin{equation}\label{EQ253}\textstyle
\sq (\P (\a ))\cong \prod _{i=0}^n \Big( \Big( \rp ^{r_i}( P(\o ^{\g _i} )/ \I _{\o ^{\g _i} })\Big)^+ \Big)^{s_i} ,
\end{equation}
where $\a=\o ^{\g _n +r_n }s_n + \dots + \o ^{ \g _0 +r_0 }s_0 +k$ is a countable ordinal presented in the Cantor normal form,
$k,r_i \in \o$, $s_i \in \N$, $\g _i \in \Lim \cup \{ 1 \}$ and $\g _n +r_n > \dots > \g _0 +r_0$.
Thus, the partial order $\sq (\P (\a ))$ is a forcing product
of iterated reduced products
of Boolean algebras $P(\o ^{\g } )/ \I _{\o ^{\g } }$,
where $\I _{\o ^{\g } }=\{ I\subset \o ^{\g }:\o ^{\g }\not\hookrightarrow I\}$ is the corresponding ordinal ideal.
Concerning the aforementioned classification, for a countable ordinal $\a \geq \o +\o$ we have $\ro (\sq (\P (\a )))\cong \ro ((P(\o )/\Fin )^+ \ast \pi)$,
where $\pi$ is a $(P(\o )/\Fin )^+$-name for a $\s$-closed separative atomless poset;
also, $\,{\mathfrak h}=\o _1$ implies that $\ro (\sq (\P (\a )))\cong \ro (P(\o )/\Fin )$, for each countable ordinal $\a$.

Unlike the countable case,
the prominent role in the classification of uncountable ordinals is played by collapsing algebras and their iterated reduced powers \cite{Kordu}.
We recall that, for cardinals $\k\geq 2$ and $\l \geq \o$,
the collapsing algebra $\Col (\l ,\k)$ is defined to be the Boolean completion  of the reversed tree $\la {}^{<\l }\k ,\supset\ra$.
So, in \cite{Kordu}, using a representation similar to (\ref{EQ253}),
it was proved that for an uncountable ordinal $\a$ we have: either the partial order $\sq (\P (\a ))$ is $\s$-closed
and completely embeds a finite power of $(P(\o )/\Fin )^+$,
or collapses at least $\o _2$ to $\o$ and, under some additional conditions, $\ro (\sq (\P (\a ) ))\cong \Col (\o ,2^{|\a|})$.

In order to present our investigation in a more general context we mention some related results.
First, taking the minimal values for $\l$ and $\k$ we obtain the collapsing algebra $\Col (\o ,2):=\ro ({}^{<\o }2)$
having several isomorphic incarnations: $\ro (\Clop (2^\o))$, $\RO (2^\o)$, $\Borel (2^\o)/\CM$.
While the important role of the algebra $\C :=\Clop (2^\o)$,
called the Cantor algebra or the Cohen algebra,
was recognized long time ago,
it turned out that the reduced power $\rp (\C )$ is a relevant structure too.
For example, using the fact that $\ro(\rp (\C ))\cong \RO (\b \R \setminus \R)$,
Dow \cite{Dow2} proved the consistency of $\RO (\b\o\setminus \o)\not\cong \RO (\b \R \setminus \R)$.
In \cite{BH} Balcar and Hru\v s\'{a}k continued the investigation of the cardinal invariants of the algebra $\rp (\C )$;
they proved that $\ft =\ft (\rp (\C)) \leq \fh (\rp (\C))\leq \min \{ \fh ,\add (\CM)\}$,
that $\ft=\fh$ implies $\RO (\b \R \setminus \R)\cong \RO (\b \o \setminus \o)$
and that $\ft <\fh (\rp (\C))$ is consistent.
The consistency of the inequality $\fh (\rp (\C))<\min \{ \fh ,\add (\CM)\}$
was established by Brendle in \cite{Bren2}.

Second, the Boolean completions of the algebras $P(\k )/[\k ]^{<\k }$, where  $\k\geq \o$ is a cardinal,
are consistently isomorphic to collapsing algebras.
In particular, under $\ft =\fh$ this is true for the algebra $P(\o)/\Fin \cong \Clop (\b \o \setminus \o)$
and for the corresponding reduced power we have $\rp (P(\o)/\Fin)\cong P(\o \times \o)/(\Fin \times \Fin)=\Clop ((\o \times \o ^\ast)^\ast)$.
Regarding the cardinal invariants of this reduced power
Szyma\'nski and Zhou Hao Xua in \cite{Szy} proved that $\ft (\rp (P(\o)/\Fin))=\o _1$ in ZFC,
Dow in \cite{Dow1} proved that $\fp =\fc$ implies that $\fh (\rp (P(\o)/\Fin))=\fc$,
while the consistency of $\fh (\rp (P(\o)/\Fin))=\o _1 <\fh =\fc$ was shown by Hern\'{a}ndez-Hern\'{a}ndez in \cite{Her}.

Our main intention is to classify the algebras $\rp ^n (\Col (\l,\k))$, $n\in \N$, up to isomorphism of their Boolean completions,
assuming that
\begin{itemize}
\item[($\ast$)] {\it $\l\geq \o$ and $\k\geq 2$ are cardinals and $\k\l>\o$.}
\end{itemize}
(For $\k\l =\o$, $\Col (\l,\k)$ is the Cantor algebra mentioned above.)
It turns out that, sometimes under additional assumptions, these completions $\ro(\rp ^n (\Col (\l,\k)))$ are collapsing algebras again.
Our strategy is to use the methods of forcing
and, as usual, to apply them to a convenient preorder which is forcing equivalent to the reduced power under consideration.
For this aim in Section \ref{S3} we detect conditions providing that the Boolean completion of a preorder is a collapsing algebra.
In Section \ref{S4} to a reduced power $\rp (\Col (\l,\k))$ we adjoin a forcing equivalent preorder $\Q _{\l,\k}$,
such that $\k ^{<\l}\leq |\Q _{\l,\k}|\leq (\k ^{<\l})^\o$,
which implies that the algebra $\rp (\Col (\l,\k))$, regarded as a forcing notion,
preserves cardinals above $(\k ^{<\l})^\o$.
In addition, it is known that it preserves $\o _1$ (since it is $\s$-closed)
and collapses $\fc$ to $\fh$ (since completely embeds $P(\o )/\Fin$).
So, regarding our intention, an important task is to detect which cardinals are collapsed and which are preserved.

It turns out that,
regarding the cardinal collapsing and the representation of the Boolean completion of the reduced power $\rp (\Col (\l ,\k))$ and its iterations,
the main role is played by the cofinality of the cardinal $\l$.
So, in Section \ref{S5} we regard the case when $\cf (\l)>\o$ and show that,
in ZFC, the algebra $\rp (\Col (\l ,\k))$ forces $|\k ^{<\l}|=|\cf (\l)|$ (as the algebra $\Col (\l ,\k)$ itself).
If $\cf (\l)=\o _1$ or, under $\fh =\o _1$, whenever $\cf (\l) \in [\o _1,\fc]$,
the Boolean completions of all iterations have a unique representation and the collapse is the maximal possible:
$\ro(\rp ^n (\Col (\l ,\k)))\cong \Col (\o _1,\k ^{<\l} )$, for all $n\in \N$.
That statement is consistently false:  $\fp =\fc >\o _1$ implies that $\fh (\rp (\Col (\fc ,\fc))=\fc$ (Dow \cite{Dow1}).
In addition, if $\fh =\o _1$, then $\ro(\rp ^n (\Col (\o_1 ,\k)))\cong \Col (\o _1,\fc )\cong \ro (P(\o)/\Fin)$,
whenever $\k \in [2,\fc]$.

Concerning the case when $\cf (\l)=\o$ in Section \ref{S6} we show that,
in ZFC, the reduced power $\rp (\Col (\l ,\k))$ forces $|\k ^{<\l}|\leq |\fh|$
and if, in addition, $\cf (\k ^{<\l})=\o$, then it forces $|(\k ^{<\l})^+|\leq |\fh|$.
The maximal collapse of $\k ^{<\l}$ or $(\k ^{<\l})^+$ (to $\o _1$) is provided by the equality $\fh =\o _1$,
but also if $\fb=\fd$ and $0^\sharp$ does not exist or if $\k ^{<\l}\geq \fc$;
then we obtain a complete picture under the Singular Cardinal Hypothesis (SCH).
In Section \ref{S7} we join our results
and, assuming that SCH and $\fh =\o _1$  hold,\footnote{
We note that the equality ${\mathfrak h}=\o _1$ follows from CH
and, moreover, it holds in many models of $\neg\,$CH;
for example, in all iterated forcing models obtained by adding Cohen, random, Sacks, Miller, Laver and Hechler reals;
it is false under MA + $\neg\,$CH but then (or, more generally, if $\fb=\fc$, which holds in the Mathias model) we have $\fb=\fd$ (see \cite{Blas}).}
show that for any cardinals $\l$ and $\k$ satisfying $(\ast)$
and $\cf (\l )\leq \fc$ and each $n\in \N$ we have
\begin{equation}\label{EQ354}
\ro (\rp ^n(\Col (\l,\k)))\cong   \left\{
                                         \begin{array}{ll}
                                              \Col (\o _1,\fc )          ,     & \mbox{ if } \k^{<\l} \leq \fc ;\\[1mm]
                                              \Col (\o _1,\k^{<\l})      ,     & \mbox{ if } \k^{<\l}    > \fc \land \cf (\k^{<\l})>\o;\\[1mm]
                                              \Col (\o _1,(\k^{<\l}) ^+ ),     & \mbox{ if } \k ^{<\l}   > \fc \land \cf (\k^{<\l})=\o.
                                         \end{array}
                                  \right.
\end{equation}
On the other hand, if  $\fb =\fd$ and $0^\sharp$ does not exist, then (\ref{EQ354}) holds when $\cf (\l )= \o$.
\section{Basic definitions and facts}\label{S2}
\noindent
{\bf Preorders and forcing.} If $\P =\la P, \leq \ra $ is a preorder, the elements $p$ and $q$ of $P$ are {\it incompatible},
we write $p\perp q$, iff there is no $r\in P$ such that $r\leq p,q$.
A set $A\subset P$ is an antichain iff $p\perp q$, for different $p,q\in A$.
An element $p$ of $P$ is an {\it atom } iff each two elements $q,r\leq p$ are compatible.
${\mathbb P} $ is called {\it atomless} iff it has no atoms.
A set $D\subset P$ is called: {\it dense} iff for each $p\in P$ there is $q\in D$ such that $q\leq p$;
{\it open} iff $q\leq p\in D$ implies $q\in D$.
If $\l $ is a cardinal, ${\mathbb P} $ is called {\it $\l$-closed} iff for each
$\gamma <\l $ each sequence $\langle  p_\alpha :\alpha <\gamma\rangle $ in $P$, such that $\alpha <\beta \Rightarrow p_{\beta}\leq p_\alpha $,
has a lower bound in $P$. $\omega _1$-closed preorders are called {\it $\sigma$-closed}.

A preorder ${\mathbb P} =\langle  P , \leq \rangle $ is called
{\it separative} iff for each $p,q\in P$ satisfying $p\not\leq q$ there is $r\leq p$ such that $r \perp q$.
The {\it separative modification} of ${\mathbb P}$
is the separative preorder $\mathop{\rm sm}\nolimits ({\mathbb P} )=\langle  P , \leq ^*\rangle $, where
$p\leq ^* q \Leftrightarrow \forall r\leq p \; \exists s \leq r \; s\leq q $; then $\sm (\sm (\P))=\sm (\P)$.
The {\it separative quotient} of ${\mathbb P}$
is the partial order $\sq (\P)=\asq (\sm (\P))$, where $\asq (\P )$ is the {\it antisymmetric quotient} of a preorder $\P$.
Then  $\ro (\sq (\P))$ is the {\it Boolean completion of $\P$}.
Preorders $\P$ and $\Q$ are called {\it forcing equivalent}, in notation $\P \equiv_{forc}\Q $,  iff they produce the same generic extensions.
\begin{fac}(Folklore)  \label{T2226}
If $\P$ is a preorder, then

(a) $\P\equiv_{forc}\asq (\P )\equiv_{forc}\sm (\P)\equiv_{forc}\sq (\P)\equiv_{forc}\ro (\sq (\P))$;

(b) $\P$ is $\l$-closed iff $\asq (\P )$ is $\l$-closed.
\end{fac}
If $\la \P ,\leq _\P , 1_\P \ra$ and  $\la \Q , \leq _\Q , 1_\Q \ra$ are preorders, then a mapping
$f: \P \rightarrow \Q$ is a {\it complete embedding}, in notation $f:\P \hookrightarrow _c \Q$,  iff

(ce1) $p_1 \leq _\P p_2 \Rightarrow f(p _1)\leq _\Q f(p_2)$,

(ce2) $p_1 \perp _\P p_2 \Leftrightarrow f(p _1)\perp _\Q f(p_2)$,

(ce3) $\forall q\in \Q \; \exists p\in \P \; \forall p' \leq _\P p \;\; f(p') \not\perp _\Q q$.

\noindent
Then $\Q\equiv_{forc}\P \ast  \pi $,
where $\pi$ is a $\P$ name for a partial order (see \cite{Kun}).
A mapping $f: \P \rightarrow \Q$ is a {\it dense embedding}, in notation $f:\P \hookrightarrow _d \Q$,
iff (ce1) and (ce2) hold and $f[P]$ is a dense suborder of $\Q$. For the following claims see \cite{Kun}.
\begin{fac}\label{T338}
If $\P_1$, $\P_2$ and $\P_3$ are separative partial orders (or preorders only), then

(a) If $f:\P_1 \hookrightarrow _c \P_2$ and  $g:\P_2 \hookrightarrow _c \P_3$, then $g\circ f:\P_1 \hookrightarrow _c \P_3$;

(b) If $f:\P_1 \hookrightarrow _d \P_2$, then $f:\P_1 \hookrightarrow _c \P_2$;

(c) If $\P_1 \hookrightarrow _c \P _2$, then $\ro (\P _1 )^+ \hookrightarrow _c \ro (\P_2 )^+$;

(d) If $\P_1 \hookrightarrow _d \P_2$, then $\ro (\P_1 ) \cong \ro (\P_2 )$.

(e) If $\P_1 ,\P_2 \hookrightarrow _c \P _1 \times \P_2$.
\end{fac}

\begin{fac}\label{T209}
(Balcar, Pelant and Simon, \cite{BPS})
Forcing by $P(\o)/\Fin$ collapses $\fc$ to $\h$.
\end{fac}
\noindent
{\bf Cardinals.} Definitions and basic facts concerning the cardinal invariants of the continuum $\fp,\ft, \fh, \fb,\fd$ and $\fu$
can be found in \cite{Blas}; we only note that  $\o _1 \leq\fp=\ft\leq \fh\leq \fb\leq \fd\leq \fc$ and that $2^\k =\fc$, whenever $\o \leq \k <\ft$.
The following fact from cardinal arithmetic will be frequently used; see \cite{Comf}, p.\ 17.
\begin{fac}\label{T327}
If $\l\geq \o$ and $\k \geq 2$ are cardinals, then

(a) $(\k ^{<\l } )^{<\cf(\l )}=\k ^{<\l }\geq \l$; so, $(\k ^{<\l } )^{<\l }=\k ^{<\l }$, if $\l$ is regular;

(b) Under {\rm GCH} we have $\k ^\o =\k $ iff $\cf (\k )>\o$.
\end{fac}
{\bf Reduced powers.}
if $\P =\la P ,\leq \ra$ is a preorder,
then the {\it reduced power} $\P^\o/\Phi$ is the quotient structure $\la P^\o/\!=_\Phi ,\leq _\Phi\ra$,
where $=_\Phi$ is the equivalence relation on the set $P^\o$ defined by $p=_\Phi q$ iff $\{ i\in \o : p_i=q_i\}\in \Phi$,
and $[p]_{=_\Phi} \leq _\Phi [q]_{=_{\Phi}}$ iff $\{ i\in \o : p_i\leq q_i\}\in \Phi$.
The reduced power of a preorder (partial order, Boolean algebra) is a preorder (partial order, Boolean algebra),
since the corresponding first-order theories are Horn theories \cite{Chang},
but the completeness of Boolean algebras is not preserved in reduced powers: the algebra $2^\o/\Phi\cong P(\o )/\Fin$ is incomplete.
For convenience, instead of the notation  $\P^\o/\Phi$ we use $\rp (\P)$
and define {\it iterated reduced powers}: $\rp ^0 (\P ):=\P $ and $\rp ^{n+1} (\P ):=\rp (\rp ^n (\P ))$, for $n\in \o$.
\begin{fac}\label{T300}
If $\B$ and $\C$ are non-trivial Boolean algebras, then

(a) $\rp(\B )^+$ is a $\s$-closed poset, $\rp(\B )^+ \equiv_{forc}(P(\o )/\Fin )^+\ast \pi$
and we have $\o _1 \leq \fh(\rp(\B )^+ )\leq  \fh $ (see \cite{BH});

(b) If $\B^+ \hookrightarrow _c \C ^+$, then $\rp (\B )^+ \hookrightarrow _c \rp (\C )^+$ (see \cite{Bren1}).
\end{fac}
\begin{fac}\label{T334}
If $\B$ and $\C$ are Boolean algebras, then we have

(a) If $\ro (\B )\cong \C$, then $\ro (\rp (\B)) \cong \ro (\rp (\C))$;

(b) $\ro (\rp ^n(\ro (\B)))\cong \ro (\rp ^n(\B))$, for each $n\in \N$.
\end{fac}
\dok
(a) Since $\B ^+$ is dense in $\ro (\B )^+\cong \C ^+$
w.l.o.g.\ we assume that $\B$ is a subalgebra of $\C$ and that $\B ^+$ is dense in $\C ^+$.
It is easy to check that the mapping, $f:\rp (\B )\rightarrow \rp (\C )$
defined by $f([\la b_i \ra]_{=_{\rp (\B )}})=[\la b_i \ra]_{=_{\rp (\C )}}$, for all $\la b_i\ra\in B^\o$,
is an embedding; thus $f[\rp (\B )]$ is a subalgebra of $\rp (\C )$.
In addition, if $[\la c_i \ra]_{=_{\rp (\C )}}\in \rp (\C ) ^+$,
then $\supp (\la c_i\ra):=\{ i\in \o : c_i > 0\}\in [\o ]^\o$,
and, since $\B ^+$ is dense in $\C ^+$, for each $i\in \supp (\la c_i\ra)$ there is $b_i\in \B ^+$ such that $b_i\leq c_i$.
Defining $b_i=0$, for $i\in \o \setminus \supp (\la c_i\ra)$,
we have $f([\la b_i \ra]_{=_{\rp (\B )}})=[\la b_i \ra]_{=_{\rp (\C )}} \leq [\la c_i \ra]_{=_{\rp (\C )}}$;
thus, the set $f[\rp (\B )^+]=f[\rp (\B )]^+$ is dense in $\rp (\C )^+$.
So, $\ro (\rp (\B))=\ro (\rp (\B )^+) \cong \ro (f[\rp (\B )^+])\cong \ro(\rp (\C )^+) =\ro (\rp (\C))$.

(b) For $n=1$ by (a) we have $\ro (\rp (\ro (\B)))\cong \ro (\rp (\B))$.
Assuming that $\ro (\rp ^n(\ro (\B)))\cong \ro (\rp ^n (\B))$ we have
\begin{eqnarray*}
\ro (\rp ^{n+1}(\ro (\B))) &    =    & \ro (\rp(\rp ^n(\ro (\B))))\quad \hfill (\mbox{by definition of }\rp^n)\\
                           &  \cong  & \ro (\rp(\ro (\rp ^n(\ro (\B)))))\quad \hfill (\mbox{by (a) applied to $\rp ^n(\ro (\B))$})\\
                           &  \cong  & \ro (\rp(\ro (\rp ^n (\B))))\quad \hfill (\mbox{by the induction hypothesis})\\
                           &  \cong  & \ro (\rp(\rp ^n (\B)))\quad \hfill (\mbox{by (a) applied to $\rp ^n (\B)$})\\
                           &  \cong  & \ro (\rp ^{n+1} (\B)))\quad \hfill \hfill (\mbox{by definition of }\rp^n) \quad \Box
\end{eqnarray*}
{\bf Collapsing algebras.}
Let $\l \geq \o$ and $\k \geq 2$ be cardinals and let ${}^{<\l }\k :=\bigcup _{\a <\l}{}^{\a }\k$.
The {\it collapsing algebra} $\Col (\l,\k)$ is the Boolean completion of the reversed tree $\la {}^{<\l }\k ,\supset \ra$.
The following fact is well known (see \cite{Balc1}, p.\ 342).
\begin{fac}\label{T368}
If $\l$ is a regular cardinal and $\k \geq 2$, then $\Col (\l ,\k)\cong \Col (\l ,\k ^{<\l})$,
the forcing $\Col (\l ,\k)^+$ collapses $\k ^{<\l }$ to $\l$,
preserves cardinals $\leq \l $ and $>\k ^{<\l}$,
contains a $\l $-closed dense suborder and $\pi (\Col (\l ,\k))=\k ^{<\l}$.
\end{fac}
The following statement (perhaps folklore)
to a collapsing algebra gives its {\it real} name, describing the whole interval of collapse (note that $(\k ^{<\l })^{<\cf(\l )}=\k ^{<\l }$).
\begin{fac}\label{T351}
If $\l\geq \o$ and $\k\geq 2$ are cardinals, then $\Col (\l ,\k)\cong \Col (\cf(\l ) ,\k ^{<\l })$.
\end{fac}
\dok
If $\l$ is a regular cardinal, the statement follows from  Fact \ref{T368}.

Otherwise, $\mu :=\cf (\l )<\l$ is a regular cardinal.
Let $\la \l _\xi : \xi <\mu \ra \in {}^{\mu }(\l \setminus \mu)$ be an increasing sequence of cardinals which is cofinal in $\l$.
Then for each $\f \in {}^{<\l}\k$ there is $\xi <\mu$ such that $\dom (\f)\subset \l _\xi$;
thus $D:=\bigcup _{\xi <\mu}{}^{\l _\xi}\k$ is a dense suborder of ${}^{<\l}\k$
and, hence, of $\Col (\l ,\k)$ as well.
By the regularity of $\mu$ the poset $D$ is $\mu$-closed
so forcing $\Col (\l ,\k)$ preserves $\mu$.
In addition, $\k ^{<\l }:=\sup \{ \k^\nu :\nu <\l\}= \bigcup _{\xi <\mu}  \k^{\l _\xi}$
and we show that $\Col (\l ,\k) \Vdash | ((\k^{\l _\xi})^V)\check{\;}|= \check{\mu}$, for each $\xi <\mu$.
For $f \in {}^{\l _\xi}\k$ let
\begin{eqnarray*}
D_f  & := & \Big\{ \p \in {}^{<\l}\k : \exists \zeta \in \mu \setminus \xi \\
      &&  \;\; \Big(\{\l _\zeta +\a :\a< \l _\zeta \} \subset \dom (\p)\land \forall \a <\l _\xi \;\; \p (\l _\zeta +\a)=f(\a ) \Big)\Big\},
\end{eqnarray*}
where $\l _\zeta +\l _\xi$ denotes the sum of ordinals.

For $\f \in {}^{<\l}\k$ there exists $\zeta \in \mu \setminus \xi$ such that $\dom (\f)\subset \l _\zeta$
and, defining $\p =\f \cup \{ \la \a ,0\ra : \a \in \l _\zeta \setminus \dom (\f)\}\cup \{ \la \l _\zeta +\a ,f(\a )\ra :\a <\l _\xi \}$,
we have $D_f \ni \p \supset \f$; so, the sets $D_f$, $f \in {}^{\l _\xi}\k$, are dense in ${}^{<\l}\k$.

Let $G$ be a ${}^{<\l}\k$ generic filter over $V$ and let $\p _f \in G\cap D_f$, for $f \in {}^{\l _\xi}\k$.
Then for each $f \in {}^{\l _\xi}\k$ there is $\zeta_f \in \mu$
such that $\{\l _{\zeta_f} +\a :\a< \l _{\zeta _f} \} \subset \dom (\p _f)$ and  $\p _f (\l _{\zeta _f} +\a)=f(\a )$, for all $\a <\l _\xi$.
Clearly, $\f \neq g$ implies that $\zeta_f \neq\zeta_g$
and, hence, $f\mapsto \zeta_f$ is an injection from  ${}^{\l _\xi}\k$ to $\mu$.
Thus $|({}^{\l _\xi}\k)^V|^{V[G]}=\mu$, for each $\xi <\mu$, which gives $|(\k ^{<\l })^V|^{V[G]}=\mu$.

Now $\mu$ is a regular cardinal,
$D:=\bigcup _{\xi <\mu}{}^{\l _\xi}\k$ is a separative atomless $\mu$-closed poset
of size $\k ^{<\l}=(\k ^{<\l})^{<\mu}$ (see Fact \ref{T327})
and collapses $\k ^{<\l}$ to $\mu$.
By Theorem \ref{T321} given in the sequel we have $\ro (D)\cong \Col (\mu ,\k ^{<\l})$
and since $D$ is dense in  $\Col (\l ,\k)$  we have $\Col (\l ,\k)\cong \Col (\mu ,\k ^{<\l })$.
\hfill $\Box$
\section{Preorders forcing equivalent to collapsing algebras}\label{S3}
In the analysis of the algebras $\rp ^n (\Col (\l,\k))$ we will work with convenient forcing equivalent preorders
and the following statement will be frequently used.
\begin{te}\label{T321}
Let $\l\geq \o$ be a regular cardinal and $\P$ a separative $\l$-closed preorder of size $\k=\k^{<\l}$.

(a) If $\k >\l$ and $1_\P\Vdash |\check{\k}|= \check{\l}$, then $\ro (\sq (\P ))\cong \Col (\l ,\k)$;

(b) If $\k =\l$ and $\P$ is atomless, then $\ro (\sq (\P ))\cong \Col (\l ,\l)$.
\end{te}
\dok
(a) First let $\P$ be a partial order, $\B:=\ro (\P)$ and w.l.o.g.\ we suppose that $\P$ is a dense suborder of $\B ^+$. We construct a dense copy of $\la {}^{<\l }\k ,\supset \ra$ in $\B ^+$.
In each extension $V_\P [G]$ we have $|G|\leq |\k |=\l$;
so, there is a $\P$-name $\tau$ such that $1_\P\Vdash ``\tau :\check{\l} \rightarrow \Gamma \mbox{ is onto}" $, where $\Gamma :=\{ \la \check{p},p\ra :p\in \P\}$ is the name for a generic filter.
The sets $D_\a :=\{ q\in \P : \exists r \in \P \; q\Vdash \la \a ,r \ra\check{\;} \in \tau  \}$, $\a\in \l$, are open  and dense in $\P$.
For $\a <\l$ and $p\in \P$ we have $p\Vdash |\check{\k}|= \check{\l}$
and, hence, there is a maximal antichain $A$ in $(\cdot ,p]_\P$ such that $A\in [D_\a ]^\k$.
Moreover, if $b\in \B^+$ and $B$ is a maximal antichain in $\P\cap (\cdot ,b]_{\B ^+}$,
then refining each $p\in B$ as above we obtain a maximal antichain $C$ in $(\cdot ,b]_{\B ^+}$ such that $C\in [D_\a ]^\k$.
Thus
\begin{equation}\label{EQ327}
\forall \a <\l\;\;\forall b\in \B^+ \;\;\exists C\in [D_\a ]^\k \;\; (C \mbox{ is a maximal antichain in }(\cdot ,b]_{\B ^+})
\end{equation}
By recursion on levels of the tree $\la {}^{<\l }\k ,\subset \ra$
we construct an (isomorphic) embedding $p:\la {}^{<\l }\k ,\supset \ra\rightarrow \B^+$ such that for each $\a <\l$
\begin{equation}\label{EQ302}
A_\a:=\{ p_\p :\p \in {}^{\a }\k \} \mbox{ is a maximal antichain in }\B^+.
\end{equation}
First let $p_\emptyset =1_\P$.
Let $0<\a\in \l$ and let $p_\f$ be constructed for all $\f \in {}^{< \a }\k$.
If $\a =\b +1$, then for $\f \in {}^{\b }\k$
by (\ref{EQ327}) there is an antichain $\{ a^{\f}_\xi : \xi <\k \}\subset D_\b$ maximal in $(\cdot ,p_{\f}]$,
we define $p_{\f  ^{\smallfrown}\xi}=a^\f_\xi$, for $\xi\in \k$,
and (\ref{EQ302}) is true.
By the definition of the set $D_\b$, for each $\f \in {}^{\b }\k$ and $\xi\in \k$ there is a unique $r_{\f  ^{\smallfrown}\xi}\in \P$ such that
\begin{equation}\label{EQ301}
p_{\f ^{\smallfrown}\xi} \Vdash \la \b , r_{\f  ^{\smallfrown}\xi}\ra\check{\;} \in \tau.
\end{equation}
If $\a$ is a limit ordinal and $\p \in {}^{\a }\k$,
then $p_{\p \upharpoonright 0}> p_{\p \upharpoonright 1}> \dots > p_{\p \upharpoonright \b}>\dots$
and the members of this sequence of successor height belong to $\P$ and form a coinitial subsequence.
Since the poset $\P$ is $\l$-closed there is $p\in \P$ such that $p\leq p_{\p \upharpoonright \b}$, for all $\b <\a$,
and, hence, $p_\p :=\bigwedge _{\b <\a}p_{\p \upharpoonright \b} \in \B ^+$.
We show that (\ref{EQ302}) is true; let $p\in \P$.
The sets $A_\b \da$, $\b <\a$, are open and dense in $\B ^+$
and, since $\P$ and, hence, $\B ^+$ is $|\a|$-distributive, the set $\bigcap_{\b <\a} A_\b \da$ is dense.
So, there is $r\in \bigcap_{\b <\a} A_\b \da$ such that $r\leq p$.
For $\b <\a$ there is a unique $\p _\b \in {}^{\b }\k$ such that $r\leq p_{\p _\b}$,
and, by the construction, there is $\p \in {}^{\a }\k$ such that $\p \upharpoonright \b =\p _\b$, for all $\b <\a$,
which implies that $r\leq p_\p =\bigwedge _{\b <\a}p_{\p \upharpoonright \b}$.
So $p\not\perp p_\p$, (\ref{EQ302}) is true, and the recursion works.

Finally we show that the set $D:=\{ p_\f :\f \in {}^{<\l }\k\}$ is dense in $\B ^+$.
Since $\P$ is dense in $\B ^+$ it is sufficient to show that for each $r\in \P$ there is $\f \in {}^{<\l }\k$ such that $p_\f\leq r$.
Let $G$ be a $\P$-generic filter over $V$ and $r\in G$.
Then $r\in \ran (\tau _G)$ and $\la \b, r\ra\in \tau _G$ for some $\b\in \l$.
So, there is $q\in G$ such that $q\Vdash \la \b ,r \ra\check{\;} \in \tau$,
which implies that $q\in D_\b$.
In addition, there is $q'\in G$ such that $q' \leq r,q$.
By (\ref{EQ302}) $A_{\b +1}:=\{ p_{\f^{\smallfrown}\xi} :\f \in {}^{\b }\k \land \xi \in \k\}$ is a maximal antichain in $\P$.
Let $\f \in {}^{\b }\k $ and $\xi <\k$ be such that $p_{\f ^{\smallfrown}\xi}\not\perp q' $
and let $s\leq p_{\f ^{\smallfrown}\xi} , q'$.
Then $s\Vdash \la \b ,r \ra\check{\;} \in \tau$
and by (\ref{EQ301}) we have $p_{\f  ^{\smallfrown}\xi} \Vdash \la \b , r_{\f  ^{\smallfrown}\xi}\ra\check{\;} \in \tau$
which implies that $s \Vdash \la \b , r_{\f  ^{\smallfrown}\xi}\ra\check{\;} \in \tau$,
and, since $1_\P\Vdash \tau :\check{\l} \rightarrow \Gamma $, we have $r=r_{\f  ^{\smallfrown}\xi}$.
Thus $p_{\f ^{\smallfrown}\xi}\Vdash \la \b ,r \ra\check{\;} \in \tau$
and, hence, $p_{\f ^{\smallfrown}\xi}\Vdash \check{r}\in \Gamma$
which, since $\P$ is separative, gives $p_{\f ^{\smallfrown}\xi}\leq r$.

If $\P$ is a preorder, then $\sm (\P )=\P$ and $\sq (\P )=\asq (\P )$.
So by (a) and (b) of Fact \ref{T2226}
$\sq (\P )$ is a $\l$-closed separative partial order and $1_{\sq (\P )}\Vdash |\check{\k}|= \check{\l}$,
which implies that $\cc (\sq (\P ))>\k$,
and, hence, $\k \leq |\sq (\P )|\leq |\P|=\k$,
which gives  $|\sq (\P )|=\k$. As above we have $\ro (\P ):=\ro (\sq (\P ))\cong \Col (\l ,\k)$.

(b) Again we suppose that $\P$ is a $\l$-closed partial order of size $\l=\l ^{<\l}$ and that $\P \subset \B ^+$ is dense, where $\B:=\ro (\P)$.
As in (a), taking a $\P$-name $\tau$ such that $1_\P\Vdash ``\tau :\check{\l} \rightarrow \Gamma \mbox{ is onto}" $,
we construct a dense copy of the reversed tree ${}^{<\l }\l$ in $\B ^+$.
The only difference is that in the proof of (\ref{EQ327}) we use the following statement.

{\it Claim.} If $\l \geq \o$ is a regular cardinal, $\P$ an atomless $\l$-closed preorder and $p\in \P$, then there is an $\l$-sized antichain in $(\cdot ,p]$.

{\it Proof of Claim.}
By recursion we define $p_\a$ and $q_\a$, for $\a <\l$, such that for each $\a,\b <\l$ we have
(i) $\a <\b \Rightarrow p_\b \leq p_\a$,
(ii) $q_{\a +1}\leq p_\a$,
(iii) $p_{\a +1}\perp q_{\a +1}$,
(iv) $\a \in \Lim _0 \Rightarrow q_\a =p_\a$.

First let $p_0=q_0 =p$. Let $0< \b <\l$ and suppose that the sequence $p _\a ,q_\a $, $\a <\b$, satisfies (i)--(iv).
If $\b =\a +1$, then, since $\P$ is atomless we can choose $p_{\a +1}, q_{\a +1} \leq p_\a$ such that $p_{\a +1}\perp q_{\a +1}$.
If $\b$ is a limit ordinal, then by (i) there is $p' \in \P$ such that $p'\leq p_\a$, for all $\a <\b$, and we define $p_\b =q_\b=p'$.
The recursion works.

Let $A:=\{ q_{\a +1}:\a <\l\}$
and suppose that there are $\a<\b <\l$ and $r\leq q_{\a +1}, q_{\b +1}$.
Then by (ii) we have $r\leq q_{\b +1} \leq p_\b$
and, since $\a +1 \leq \b$, by (i) we have $p _\b \leq p _{\a +1}$;
thus, $r\leq p _{\a +1},q _{\a +1}$, which is false by (iii).
Thus $A\subset (\cdot ,p]$ is an antichain and $|A|=\l$.
\hfill $\Box$
\begin{cor}\label{T365}
If $\P$ is a separative atomless $\s$-closed preorder of size $\mu=\mu^{\o}\geq \o _1$ and $1_\P\Vdash |\check{\mu }|= \check{\o _1}$, then $\ro (\sq (\P ))\cong \Col (\o _1 ,\mu)$.
\end{cor}
\begin{fac}\label{T367}
If $\ft =\fh$, then $\ro (P(\o)/\Fin)\cong\Col (\fh ,\k)$, for each cardinal $\k \in [2,\fc]$.
\end{fac}
\dok
Since $\ft$ is a regular cardinal and $2^\mu =\fc$, for all $\mu \in [\o ,\ft )$ (see \cite{Blas}),
we have  $\k ^{<\ft}=\fc$ and, by Fact \ref{T351}, $\Col (\fh ,\k) \cong\Col(\fh ,\fc )$.
The poset $(P(\o)/\Fin)^+ $ is separative, atomless, $\h$-closed and of size $\fc$,
by Fact \ref{T209} we have $(P(\o)/\Fin)^+ \Vdash |\check{\fc }|=\check{\fh }$
and, by Theorem \ref{T321}, $\ro (P(\o)/\Fin)\cong \Col (\fh ,\fc)$ indeed.
\hfill $\Box$
\section{The preorder {\boldmath$\Q _{\l ,\k}$}}\label{S4}
We recall that we investigate the reduced powers $\rp (\Col (\l ,\k))$ assuming ($\ast$):
$\l\geq \o$ and $\k\geq 2$ are cardinals and $\k\l>\o$.
W.l.o.g.\ we assume that ${}^{<\l }\k \subset \Col (\l ,\k)$
and by $s$ denote the zero  of $\Col (\l ,\k)$.
For $\p= \la \p _i :i\in \o \ra \in (\{ s\}\cup {}^{<\l }\k )^\o$ let us define $\supp (\p ):=\{ i\in \o : \p _i > s\}$
and let $\Q _{\l ,\k}:=\la Q, \leq _Q\ra$, where
\begin{eqnarray}
Q                 & :=              & \big\{ \p \in (\{ s\}\cup {}^{<\l }\k )^\o :\;|\supp (\p )|=\o\big\},\label{EQ242}\\
\p \leq _Q \eta  & \Leftrightarrow & \exists j\in \o \;\;\forall i\in \supp (\p )\setminus j \;\;\big(i\in \supp (\eta ) \land \p _i \supset \eta _i\big).\label{EQ243}
\end{eqnarray}
\begin{prop}\label{T301}
Under $(\ast)$ for the preorder $\Q _{\l ,\k}$ we have

(a) $\Q _{\l ,\k}$ is separative atomless $\sigma$-closed and $\k ^{<\l}\leq |\Q _{\l ,\k}|\leq (\k ^{<\l})^\o$;

(b) $ \ro (\rp (\Col (\l ,\k)))\cong\ro (\sq (\Q _{\l ,\k}))$ and, hence, $\rp (\Col (\l ,\k))^+\equiv_{forc}\Q _{\l ,\k}$.
\end{prop}
\dok
(a) Clearly, $\Q _{\l ,\k}$ is atomless; for a proof that it is separative we take $\p,\eta \in Q$, where $\p \not\leq _Q \eta$.
By (\ref{EQ242}) we have $\supp (\p ),\supp (\eta )\in [\o ]^\o$
and, by (\ref{EQ243}), $S_1 \cup S_2\in [\o ]^\o$,
where $S_1:=\supp (\p )\setminus \supp (\eta )$
and $S_2:=\{ i\in \supp (\p ) : \p _i \not \supset \eta _i\}$.

If $S_1\in [\o ]^\o$,
let $\p' _i =\p _i$, for $i\in S_1$; and $\p' _i =s$, for $i\in \o\setminus S_1$.
Then $\p' \leq _Q \p$
and, assuming that there is $\theta \in Q$ such that $\theta \leq _Q \p',\eta$,
by (\ref{EQ243}) we would have $[\o ]^\o \ni\supp (\theta)\subset ^*\supp (\p' )\cap \supp (\eta )=\emptyset$,
which is false. Thus $\p'\perp _Q \eta$.

If $|S_1|<\o$, then we have $\supp (\p )\subset^* \supp (\eta )$, $|S_2|=\o$,
and, hence, $S_2':=\{ i\in \supp (\p )\cap  \supp (\eta ): \p _i \not \supset \eta _i\}\in [\o ]^\o$.
Since the reversed tree ${}^{<\l }\k$ is separative,
for $i\in S_2'$ we take $\p _i ' \supset \p _i$ such that $\p_i' \perp \eta _i$,
and define $\p' _i =s$, for $i\in \o\setminus S_2'$.
Then $Q\ni\p' \leq _Q \p$
and, assuming that $\theta \leq _Q \p',\eta$, for some $\theta \in Q$,
by (\ref{EQ243}) we would have $[\o ]^\o \ni\supp (\theta)\subset ^* S_2' $,
and there would be $i\in S_2'$ such that $\theta _i \supset \p _i', \eta _i$,
which is false. Thus $\p'\perp _Q \eta$ again and $\Q _{\l ,\k}$ is separative indeed.

We prove that $\Q _{\l ,\k}$ is $\s$-closed.
If $\p^n=\la \p^n _i : i\in \o\ra \in Q$, for $n <\o$,
and for each $n\in \o$ we have $\p^{n+1}\leq _Q \p^n$, then by (\ref{EQ243}) there is $j_n\in \o$ such that
\begin{equation}\label{EQ306}
\forall i\in \supp (\p ^{n+1} )\setminus j_n \;\;(i\in \supp (\p ^n ) \land \p^{n+1} _i \supset \p^n _i).
\end{equation}
By recursion we define $i_n\in \o$ and $\eta _{i_n}\in {}^{<\l}\k$, for $n\in \o$, as follows. Let $i_0 =\min (\supp (\p ^1) \setminus j_0)$ and $\eta _{i_0}=\p ^1_{i_0}$. For $n\geq 1$ let
\begin{equation}\label{EQ307}
i_n :=\min \big\{ i\in \supp (\p ^{n+1}): i>i_{n-1}, j_0,\dots,j_n \big\}  \quad\mbox{ and }\quad \eta _{i_n}=\p ^{n+1}_{i_n}.
\end{equation}
Defining $S:=\{ i_n :n\in \o\}$, $\eta _i=s$, for $i\in \o \setminus S$, and $\eta =\la \eta _i :i\in \o\ra$,
by (\ref{EQ307}) we have $S=\supp (\eta)\in [\o ]^\o$ and, hence, $\eta \in Q$.

By (\ref{EQ307}) we have $i_n\in \supp (\p ^{n+1} )\setminus j_n$ and by (\ref{EQ306}) $i_n\in \supp (\p ^n )$ and $\p^{n+1} _{i_n} \supset \p^n _{i_n}$;
so, by (\ref{EQ307}) again, $i_n\in \supp (\p ^{n} )\setminus j_{n-1}$ and, by (\ref{EQ306}),  $i_n\in \supp (\p ^{n-1} )$ and $\p^{n} _{i_n} \supset \p^{n-1} _{i_n}$.
Continuing in this way we obtain $i_n\in \supp (\p ^k)$, for all $k\leq n+1$, and $\eta _{i_n}=\p ^{n+1}_{i_n}\supset\p ^n_{i_n} \supset \dots \supset \p ^0_{i_n}$.
In other words, for each $k>0$ and $n\geq k-1$ we have $i_n\in \supp (\p ^k)$ and $\eta _{i_n}\supset\p ^k_{i_n}$;
so, for each $i\in \supp (\eta)\setminus \{ i_n : n<k-1\}$ we have $i\in \supp (\p ^k)$ and $\eta _{i}\supset\p ^k_{i}$,
which by (\ref{EQ243}) means that $\eta \leq _Q \p ^k$.
Thus, $\eta$ is a lower bound for the sequence $\la \p^n :n \in\o\ra$.

Taking the constant sequences $\la \f\ra$ with support $\o$, for $\f\in {}^{<\l}\k$ we have  $|\Q _{\l ,\k}|\geq \k ^{<\l}$;
finally, $|\Q _{\l ,\k}|\leq |(\{s\}\cup {}^{<\l}\k)^\o|= (\k ^{<\l})^\o$.

(b) Clearly, $\B :=\rp (\Col (\l ,\k))$ is a Boolean algebra; for $b=\la b_i:i\in \o\ra\in \Col (\l ,\k)^\o$ let $\supp (b):=\{ i\in \o :b_i>s\}$.
Then we have $0_\B =[\la s,s,\dots \ra]_\B =\{ b\in \Col (\l ,\k)^\o: |\supp (b)|<\o \}$.
Regarding $\{ s\}\cup {}^{<\l }\k$ as a suborder of $\Col (\l,\k)$,
the reduced power $\BT := \rp (\{ s\}\cup {}^{<\l }\k)$ is a partial order with a smallest element
$0_\BT =[\la s,s,\dots \ra]_\BT=\{ \p\in (\{ s\}\cup {}^{<\l }\k)^\o: |\supp (\p)|<\o\}$.

Defining $\B ^+ :=\B \setminus \{0_\B \}$ and  $\BT ^+ :=\BT \setminus \{0_\BT \}$
we show first that the mapping $f: \BT^+ \rightarrow \B ^+$
given by $f([\p ]_\BT)=[\p ]_\B$, for all $\p\in(\{ s\}\cup {}^{<\l }\k)^\o$,
is an  isomorphism onto a dense suborder of $\B ^+$,
which will give $\ro(\BT ^+ )\cong \ro( \rp (\Col (\l ,\k)) )$.
Let $[\p]_\BT,[\p']_\BT\in \BT ^+$.
Then $[\p]_\BT=[\p']_\BT$
iff $\{ i\in\o : \p_i =\p_i'\}\in \Phi$, where $\Phi$ is the Fr\'{e}chet filter in $P(\o)$,
iff $[\p]_\B =[\p ']_\B $;
thus, the mapping $f$ is well defined and one-to-one.
In addition, $|\supp (\p)|=\o$,
which implies that $[\p]_\B \neq 0_\B $
so $f$ maps $\BT ^+$ into $\B ^+$ indeed.
Also, $[\p]_\BT \leq [\p ']_\BT$
iff $\{ i\in\o : \p_i \leq \p_i'\}\in \Phi$
iff $[\p ]_\B \leq [\p ']_\B$
iff $f([\p]_\BT )\leq f([\p ']_\BT )$;
thus, $f$ is an embedding of posets.
Further we show that the set $f[\BT ^+] $ is dense in $\B ^+$.
If $[b]_\B \in \B ^+$,
then $\supp (b)\in [\o ]^\o$
and, since  the set ${}^{<\l }\k$ is dense in $\Col (\l ,\k)^+$,
for each $i\in \supp (b)$ we can choose $\p_i\in {}^{<\l }\k$ such that $\p_i\leq b_i$.
Defining $\p_i =s$, for  $i\in \o \setminus\supp (b)$,
we obtain $[\p]_\BT\in\BT ^+$
and $f([\p]_\BT)=[\p]_\B \leq [b]_\B $.

Thus $\ro(\BT ^+)\cong \ro(\rp (\Col (\l ,\k)))$ and we show that $\BT ^+ =\sq (\Q _{\l ,\k})$.
Since the preorder $\Q _{\l ,\k}$ is separative we have $\sq (\Q _{\l ,\k})=\asq (\Q _{\l ,\k})$.
By (\ref{EQ242}) we have $\BT ^+=\{ [\p ]_\BT :\p\in (\{ s\}\cup {}^{<\l }\k)^\o \land |\supp (\p)|=\o \}=\{ [\p ]_\BT :\p\in Q \}$
and we first show that for $\p ,\eta \in Q$ we have
$\p \leq _Q \eta$ iff $[\p ]_\BT \leq [\eta ]_\BT$ (that is, iff $S:=\{ i\in \o : \p _i \leq \eta _i \}\in \Phi$).
So, if $\p \leq _Q \eta$ and $j\in\o$ satisfies (\ref{EQ243}),
then for $i\in \o \setminus j$ we have: $\p _i \supset \eta _i$, if $i\in \supp (\p)$;
and $\p _i =s \leq \eta _i$, if $i\not\in \supp (\p)$;
thus, $\o \setminus j\subset S\in \Phi$.
Conversely, if $S\in \Phi$,
then there is $j\in \o$ such that $\p _i \leq \eta _i $, for all $i>j$, and we prove (\ref{EQ243}).
So, if $i\in \supp (\p )\setminus j$, then $s<\p _i \leq \eta _i$,
which implies that $i\in \supp (\eta )$ and $\p _i \supset \eta _i$,
and the equivalence is proved.
Now, since $\BT$ is a partial order, for $\p ,\eta \in Q$ we have
$[\p ]_{\asq (\Q)}=[\eta ]_{\asq (\Q)}$ iff  $\p \leq _Q \eta \land \eta \leq _Q \p$  iff  $[\p ]_\BT = [\eta ]_\BT$.
Also, since two quotient structures are well defined,
we have $[\p ]_{\asq (\Q)} \leq_{\asq (\Q)} [\eta ]_{\asq (\Q)}$ iff  $\p \leq _Q \eta $  iff  $[\p ]_\BT \leq [\eta ]_\BT$.
\hfill $\Box$
\section{{\boldmath$\l $}'s of uncountable cofinality}\label{S5}
\subsection{Collapse of cardinals}
\begin{te}\label{T352}
$(\ast)$ If $\;\cf (\l)> \o $, then

(a) $\rp(\Col (\l,\k))^+\Vdash |((\k ^{<\l })^V)\check{\;}|=|(\cf ^V (\l))\check{\;}|$;

(b) $\rp(\Col (\l,\k))^+\!\Vdash \!|((\k ^{<\l })^V)\check{\;}|=\check{\o _1}$, if $\cf (\l) \!\leq \!\fc$ and $\min \{ \fh , \cf (\l)\}=\o _1$;\footnote{
More generally, if $\rp(\Col (\l,\k))^+\Vdash |(\cf ^V (\l))\check{\;}|\leq|(\fc ^V)\check{\;}|\land \min \{|(\fh ^V)\check{\;}|, |(\cf ^V (\l))\check{\;}|\}=\check{\o _1}$.}

(c) $\o _1 \leq \fh(\rp (\Col (\l ,\k))^+ )\leq \min \{ \fh , \cf (\l ) \}$.
\end{te}
\dok
(a) Let $\l_0 := \cf (\l)\geq \o _1$ and $\k _0:=\k^{<\l }$.
By Fact \ref{T351} we have $\Col (\l,\k)\cong \Col (\l _0 ,\k_0)$;
so we prove that $\rp(\Col (\l_0,\k_0))\Vdash |\check{\k _0}|=|\check{\l _0}|$.
Again we work with the preorder $\Q _{\l _0 ,\k _0}$
and first we show that for each $\xi <\l_0$ the set
$$\textstyle
D_\xi :=\Big\{ \p \in \Q _{\l_0 ,\k_0} : \supp (\p )=\o \land \xi \subset \bigcap _{i\in \o}\dom (\p _i)\Big\}
$$
is pre-dense in $\Q _{\l_0 ,\k_0}$.
Namely, for $\eta \in \Q _{\l_0 ,\k_0}$ we construct $\p \in D_\xi$ compatible with $\eta$.
For $i\in \supp (\eta)$ we have $\dom (\eta _i)\in \l_0$;
so, if $\xi \subset \dom (\eta _i)$, we define $\p _i =\eta _i$;
otherwise, for $\p _i$ we take an extension of $\eta _i$ to $\xi$.
For $i\in \o \setminus \supp (\eta)$ we take arbitrary $\p _i:\xi \rightarrow \k_0$;
so we have $\p\in D_\xi $.
Let $\f \in \Q _{\l_0 ,\k_0}$, where $\supp (\f)=\supp (\eta)$ and $\f _i=\p _i$, for $i\in \supp (\eta)$.
Then $\f \leq \eta ,\p$ and $\eta\not\perp \p$.

Let $G$ be a $\Q _{\l_0 ,\k_0}$-generic filter over $V$.
In $V[G]$ we take $\p ^\xi \in G \cap D_\xi$, for $\xi <\l_0$,
and define
$$\textstyle
A:=\bigcup _{\xi <\l_0} \bigcup _{i\in \o }\ran (\p ^\xi _i)\subset \k_0.
$$
For each $\xi <\l_0$ and $i\in \o $ we have $|\ran (\p ^\xi _i)|^{V[G]}\leq |\dom (\p ^\xi _i)|^{V[G]}\leq|\l_0 |^{V[G]}$;
thus, $A$ is the union of $|\l_0 \o|^{V[G]}=|\l_0|^{V[G]}$-many sets of size $\leq |\l_0|^{V[G]}$
and, hence, $|A|^{V[G]}\leq |\l_0|^{V[G]}$.

It is evident that the sets $\Delta _\zeta :=\{ \eta \in \Q _{\l_0 ,\k_0} : \forall i\in \supp (\eta )\;\; \zeta \in \ran (\eta _i) \}$, for $\zeta <\k_0$, are dense in $\Q _{\l_0 ,\k_0}$.
Let $\zeta <\k_0$ and $\eta \in G\cap \Delta _\zeta$.
Then $\dom (\eta _i)<\l_0$, for all $i\in \supp (\eta )$,
so, since $\eta\in V$ and in $V$ we have $\l _0=\cf(\l_0 )>\o$,
there is $\xi < \l_0$ such that
\begin{equation}\label{EQ325}
\forall i\in \supp (\eta ) \;\;\dom (\eta _i) \subset \xi.
\end{equation}
Since $\eta, \p ^\xi \in G$ there is $\theta \in G$ such that $\theta \leq \eta, \p ^\xi $
and, hence, there is $j\in \o$ such that
\begin{equation}\label{EQ326}
\forall i\in \supp (\theta )\setminus j \;\;(i\in \supp (\eta ) \land \theta _i \supset \eta _i ,\p ^\xi _i)
\end{equation}
(recall that $\supp (\p ^\xi)=\o$).
Now for $i\in \supp (\theta )\setminus j$,
by (\ref{EQ325}), (\ref{EQ326}) and since $\p ^\xi \in D_\xi$, we have $\dom (\eta _i) \subset \xi \subset \dom (\p^\xi _i)\subset \dom (\theta _i)$.
So by (\ref{EQ326}) we have $\eta _i \subset\p ^\xi _i$,
which, since $\eta \in \Delta _\zeta$, gives $\zeta \in \ran (\eta _i)\subset \ran (\p^\xi _i)\subset A$.

Thus $A=\k_0$,
which gives $|\k_0|^{V[G]}=|A|^{V[G]}\leq |\l_0|^{V[G]}$.
In $V$, since $\k \geq 2$, by Fact \ref{T327} we have $\k _0:=\k ^{<\l}\geq \l \geq \l _0$,
thus $|\k_0|^{V[G]}\geq |\l_0|^{V[G]}$
and the equality $|\k_0|^{V[G]}= |\l_0|^{V[G]}$ is proved.

(b) follows from (a) and Facts \ref{T209} and \ref{T300}(a).

(c) By Fact \ref{T300}(a) we have $\o _1 \leq \fh(\rp (\Col (\l ,\k))^+)\leq  \fh $.
If $\cf (\l)<\k ^{<\l}$, then by (a) forcing by $\rp (\Col (\l ,\k))^+$ adds a surjection $f:\cf (\l)\rightarrow \k ^{<\l}$
and, clearly, $f\not \in V$; thus, $\fh(\rp (\Col (\l ,\k))^+ )\leq \cf (\l )$.
Otherwise by Fact \ref{T327}(a) we have $\k ^{<\l}\geq \l \geq \cf (\l)$
and, hence, $\k ^{<\l}=\l= \cf (\l)> \o $.
By Proposition \ref{T301}(a),  $\k ^{<\l} \leq |Q_{\l ,\k}|\leq (\k ^{<\l})^\o \leq (\k ^{<\l})^{<\l }=\k ^{<\l}$,
which gives $|Q_{\l ,\k}|= \k ^{<\l}$.
So, by Proposition \ref{T301}(b), in each generic extension by $\Q_{\l ,\k}$ or, equivalently, by $\rp (\Col (\l ,\k))^+$
the set $Q_{\l ,\k}$ obtains a new subset--the $\Q_{\l ,\k}$-generic filter
and, since $|Q_{\l ,\k}|=\cf (\l)$, we have $\fh(\rp (\Col (\l ,\k))^+ )\leq \cf (\l )$ again.
\kdok
If $(\ast)$ holds and $\cf (\l)\in [\o _1,\fc]$,
then by (b) of Theorem \ref{T352} the equality $\min \{ \fh , \cf (\l)\}=\o _1$
provides a maximal collapse in extensions by  $\rp(\Col (\l,\k))$:
the cardinal $\k ^{<\l}$ is collapsed to $\o _1$.
If $\mu :=\min \{ \fh , \cf (\l)\}>\o _1$, then $\k ^{<\l}$ is collapsed at least to $\mu$
and it is natural to ask what is going on with $\mu$.
The following example shows that it is consistent that $\mu =\k ^{<\l}>\o _1$ is not collapsed at all.
\begin{ex}\label{EX300}\rm
Dow \cite{Dow1} proved that under $\fp =\fc$ we have $\fh (\rp (P(\o)/\Fin))=\fc$.
So, if $\fp =\fc >\o _1$, then $\ft =\fh=\fc$,
by Fact \ref{T367} we have $\rp (P(\o)/\Fin) \cong \rp(\Col (\fc ,\fc))$
and, hence, $\fh (\rp(\Col (\fc ,\fc)))=\fc $,
which implies that $\rp(\Col (\fc ,\fc))$ does not collapse $\min \{ \fh , \cf (\l)\}= \fc $
($\fc$ is regular since $\t$ is; see \cite{Blas}).
We note that $\ft (\rp (P(\o)/\Fin))=\o _1$ holds in ZFC  \cite{Szy}
so, in our case, $\ft (\rp (\Col (\fc ,\fc)))=\o _1$.
Moreover, by \cite{Dow1}, if $T \subset \rp (\Col (\fc ,\fc))$ is a dense tree, then $T$ is not $\o _2$-closed.
\end{ex}
\subsection{Boolean completions}
\begin{te}\label{T359}
($\ast$) Let $\;\rp (\Col (\l ,\k))^+\Vdash |((\k ^{<\l })^V)\check{\;}|=\check{\o _1}$ and
\begin{equation}\label{EQ331}
(\l =\o \mbox{ and } \k ^\o =\k)\mbox{ or }(\l >\o \mbox{ and }\cf (\l)\cf(\k ^{<\l})>\o ).
\end{equation}
Then for each $n\in \N$ we have
\begin{equation}\label{EQ332}
\ro (\rp ^n (\Col (\l ,\k)))\cong \Col (\o _1 ,\k ^{<\l}).
\end{equation}
\end{te}
\dok
Let $\f (\l,\k)$ denote the conjunction of all assumptions on $\l$ and $\k$ (i.e.\ $(\ast)$, the collapse and (\ref{EQ331})). By induction we prove that for each $n\in \N$ we have
\begin{equation}\label{EQ333}
\forall \l,\k\in \Card \;\;\Big(\f (\l ,\k)\Rightarrow \ro (\rp ^n (\Col (\l ,\k)))\cong \Col (\o _1 ,\k ^{<\l})\Big).
\end{equation}
First let $n=1$ and let $\l$ and $\k$ be cardinals satisfying  $\f (\l ,\k)$.
By Fact \ref{T351} and Proposition \ref{T301}(b) we have $\ro (\rp (\Col (\l ,\k)))\cong\ro (\rp (\Col (\cf(\l ),\k ^{<\l})))\cong \ro (\sq (\Q _{\cf(\l ),\k ^{<\l}}))$;
so, by Fact \ref{T2226}(a) we have $\rp (\Col (\l ,\k))^+ \equiv_{forc}\Q _{\cf(\l ),\k ^{<\l}}$ and, hence, $\Q _{\cf (\l), \k ^{<\l }}\Vdash |((\k ^{<\l })^V)\check{\;}|=\check{\o _1}$.
For a proof that $\ro (\sq (\Q _{\cf(\l ),\k ^{<\l}}))\cong \Col (\o _1 ,\k ^{<\l})$ we use Corollary \ref{T365},
where we replace  $\mu$ by $\k ^{<\l}$ and $\P$ by $\Q _{\cf(\l ),\k ^{<\l}}$.
So, by Proposition \ref{T301}(a) and since $\Q _{\cf (\l), \k ^{<\l }}\Vdash |((\k ^{<\l })^V)\check{\;}|= \check{\o _1}$ and by ($\ast$) we have $\k ^{<\l} \geq \o _1$,
it remains to be proved that $|\Q _{\cf(\l ),\k ^{<\l}}|=\k ^{<\l}=(\k ^{<\l})^\o $.
By Fact \ref{T327} and Proposition \ref{T301}(a) we have $\k ^{<\l}=(\k ^{<\l})^{<\cf (\l)}\leq|\Q _{\cf(\l ),\k ^{<\l}}|\leq ((\k ^{<\l})^{<\cf (\l)})^\o=(\k ^{<\l})^\o$;
so, we show that $(\k ^{<\l})^\o = \k ^{<\l}$.
If $\l =\o$ and $\k ^\o =\k$ this is evident.
If $\cf (\l )>\o$, by Fact \ref{T327} we have $(\k ^{<\l})^\o \leq (\k ^{<\l})^{<\cf (\l)}= \k ^{<\l}$ and we are done.
Finally, let $\l>\cf (\l )=\o$ and $\cf (\k ^{<\l})>\o $
and let $\la \l _n :n\in \o\ra$ be an increasing sequence of infinite cardinals cofinal in $\l$.
Since $\k ^{<\l}=\sup \{ \k ^{\l _n}:n\in \o\}$ there is $n\in \o$ such that $\k ^{<\l}=\k ^{\l _n}$
(otherwise we would have $\cf (\k^{<\l })=\o$).
Thus $(\k ^{<\l})^\o=\k ^{\l _n \o}=\k ^{\l _n}=\k ^{<\l}$ and (\ref{EQ333}) is true for $n=1$.

Assuming that (\ref{EQ333}) is true (for $n$) we prove that it is true for $n+1$.
Let $\l$ and $\k$ be cardinals satisfying  $\f (\l ,\k)$.
Then by our assumption we have  (\ref{EQ332})
and first prove that
\begin{equation}\label{EQ334}
\ro( \rp (\Col (\o _1, \k ^{<\l})))\cong\Col (\o _1, \k ^{<\l}).
\end{equation}
We have proved that (\ref{EQ333}) is true for $n=1$.
In addition, $\f (\o _1,\k^{<\l})$ is true,
because by Theorem \ref{T352}(a) we have  $\rp(\Col (\o _1,\k^{<\l}))^+\Vdash |(((\k^{<\l}) ^{<\o _1 })^V)\check{\;}|=\check{\o _1}$
and $\o _1=\cf (\o _1)>\o$ (so the second part of (\ref{EQ331}) is true).
So, by (\ref{EQ333}) (for $n=1$, replacing $\l$ by $\o _1$ and $\k$ by $\k ^{<\l}$) we obtain
$\ro (\rp (\Col (\o _1 ,\k^{<\l})))\cong \Col (\o _1 ,(\k^{<\l}) ^{<\o _1})$.
In the proof of (\ref{EQ333}) for $n=1$ we have shown that $\f (\l ,\k)$ implies $(\k ^{<\l})^\o = \k ^{<\l}$.
So, since $\f (\l ,\k)$ holds and $(\k^{<\l}) ^{<\o _1 }=(\k^{<\l}) ^\o=\k^{<\l}$, we have (\ref{EQ334}).
Finally,
\begin{eqnarray*}
\ro(\rp ^{n+1}(\Col (\l,\k))) & =     & \ro( \rp (\rp ^n(\Col (\l,\k))))        \quad \hfill (\mbox{by definition of }\rp^n)\\
                              & \cong & \ro( \rp (\ro (\rp ^n(\Col (\l,\k)))))  \quad \hfill (\mbox{by Fact \ref{T334}(a)})\\
                              & \cong & \ro( \rp (\Col (\o _1, \k ^{<\l})))     \quad \hfill (\mbox{by (\ref{EQ332})})\\
                              & \cong & \Col (\o _1, \k ^{<\l})                 \quad \hfill (\mbox{by (\ref{EQ334})}).  \quad \Box
\end{eqnarray*}
\begin{te}\label{T340}
Let $(\ast)$ hold.

(a) If $\cf (\l)=\o _1$ or ($\fh =\o _1$ and $\cf (\l) \in [\o _1,\fc]$),\footnote{More generally,
if  $\rp(\Col (\l,\k))^+\Vdash |(\cf ^V (\l))\check{\;}|\leq|(\fc ^V)\check{\;}|\land \min \{|(\fh ^V)\check{\;}|, |(\cf ^V (\l))\check{\;}|\}=\check{\o _1}$
and $\cf(\l )>\o$.}\\
then $\ro (\rp ^n (\Col (\l ,\k)))\cong \Col (\o _1 ,\k ^{<\l})$, for all $n\in \N$;

(b) If $\fh =\o _1 $ and $\k \leq \fc$, \\
then $\ro(\rp ^n (\Col (\o _1,\k)))\cong\Col (\o _1 ,\fc)\cong \ro (P(\o)/\Fin)\cong\ro(\rp ^n (P(\o)/\Fin))$.
\end{te}
\dok
(a) Under the assumptions (\ref{EQ331}) from Theorem \ref{T359} is true,
by Theorem \ref{T352}(b) we have $\rp(\Col (\l,\k))^+\Vdash |((\k ^{<\l })^V)\check{\;}|=\check{\o _1}$ and we apply Theorem \ref{T359}.

(b) Since $\k \leq \fc$ we have $\k ^\o =\fc$; so, by (a), $\ro(\rp ^n (\Col (\o _1,\k)))\cong \Col (\o _1, \fc)$.
By Fact \ref{T367} $\Col (\o _1 ,\fc)\cong \ro (P(\o)/\Fin)\cong \Col (\o _1,\k)$
and, by Fact \ref{T334}(b), $\ro(\rp ^n (\Col (\o _1,\k)))\cong \ro(\rp ^n (\ro (P(\o)/\Fin)))\cong \ro(\rp ^n (P(\o)/\Fin))$.
\hfill $\Box$
\section{{\boldmath$\l $}'s of countable cofinality}\label{S6}
\subsection{Collapse of cardinals}
\begin{te}\label{T233}
$(\ast)$ If $\cf (\l)=\o$, then

(a) $\rp(\Col (\l,\k))^+ \Vdash |((\k ^{<\l })^V)\check{\; }|\leq |(\fh ^V )\check{\;}|$;

(b) $\rp(\Col (\l,\k))^+ \Vdash |((\k ^{<\l })^V)\check{\; }|= \check{\o} _1$,
if $\rp(\Col (\l,\k))^+ \Vdash |(\fh ^V )\check{\;}|=\check{\o}_1$.
\end{te}
\dok
First, we note that $\k _0:=\k^{<\l }>\o$;
namely, if $\l>\o$, then by Fact \ref{T327} we have $\k _0 >\o$,
and if $\k >\o$, then, clearly, $\k _0 >\o$ again.
By Fact \ref{T351} we have $\Col (\l,\k)\cong \Col (\o ,\k_0)$;
so, we prove that $\rp(\Col (\o,\k_0))\Vdash |\check{\k _0}|\leq |(\fh ^V )\check{\;}| $
and we work with the preorder $\Q _{\o ,\k _0}$ from Proposition \ref{T301}.
First we show that $\rp(\Col (\o,\k_0))\Vdash |\check{\k _0}|\leq |(\fc ^V )\check{\;}| $.

In the ground model $V$,
the set $\CS :=\bigcup _{S\in [\o ]^\o} \{ S\}\times {}^{S}\o $ is the union of $\fc$-many sets of size $\fc$,
and, hence, $|\CS|=\fc$.
For $\p \in \Q _{\o ,\k _0}$ we have $S:=\supp(\p)\in [\o ]^\o$
and $\la \dom (\p  _i) :i\in \supp (\p)\ra \in {}^{S}\o$;
thus the function $f: \Q _{\o ,\k _0}\rightarrow \CS$
given by $f(\p )=\la \supp (\p), \la \dom (\p _i) :i\in \supp (\p)\ra\ra$, for $\p \in \Q _{\o ,\k _0}$, is well defined.

The sets $\Delta _\xi :=\{ \p\in Q _{\o ,\k _0} : \forall i\in \supp (\p ) \;\; \exists \eta _i \in {}^{<\o }\k_0\;\; \p _i =\eta _i{}^{\smallfrown}\xi\}$, $\xi <\k_0$, are dense in $\Q _{\o ,\k _0}$.
Namely, if $\eta\in Q _{\o ,\k _0}$, defining $\p _i=\eta _i{}^{\smallfrown}\xi$, for $i\in \supp (\eta )$, and $\p _i=s$, otherwise,
we have $\Delta_\xi\ni\la \p _i :i\in \o \ra\leq  \eta$.

Let $G$ be a $\Q _{\o ,\k _0}$-generic filter
and, in $V[G]$, let $\p^\xi \in G\cap \Delta _\xi$, for $\xi <\k _0$.
Let $g:\k_0 \rightarrow \CS$, where $g(\xi )=f(\p^\xi)=\la \supp (\p^\xi), \la \dom (\p ^\xi_i) :i\in \supp (\p^\xi)\ra\ra$, for $\xi <\k_0$.
We show that $g$ is an injection, which gives $|\k_0|^{V[G]}\leq |\CS|^{V[G]}=|\fc ^V|^{V[G]}$.
Otherwise, there would be $\xi <\zeta <\k_0 $ and $\la S,\la n_i +1 :i\in S\ra\ra \in \CS$
such that $F(\xi )=F(\zeta )=\la S,\la n_i +1 :i\in S\ra\ra $,
that is, $\supp (\p^\xi) =\supp (\p^\zeta)=S$ and $\dom (\p ^\xi_i ) = \dom (\p ^\zeta _i )=n_i +1 \in \o $, for $i\in S$.
So, for $i\in S$ we would have $\p ^\xi_i =\eta _i {}^{\smallfrown}\xi$ and $\p ^\zeta_i=\theta _i {}^{\smallfrown} \zeta$,
where $\eta _i =\p ^\xi_i \upharpoonright n_i $ and $\theta _i =\p ^\zeta_i \upharpoonright n_i$.
Since $\p^\xi , \p^\zeta \in G$ there is $\p\in G$ such that $\p\leq   \p^\xi , \p ^\zeta$
and by (\ref{EQ243}) there would be $j\in \o$
such that for each $i\in S\setminus j$ we have $\p _i \supset \p^\xi_i , \p ^\zeta _i$,
and, hence, $\xi =\p^\xi_i (n_i)=\p_i (n_i)=\p^\zeta_i (n_i)=\zeta$,
which is false.
Thus, $|(\k^{<\l })^V|^{V[G]}=|\k_0|^{V[G]}\leq |\fc ^V|^{V[G]}$.

By Fact \ref{T300}(a) we have $\rp(\Col (\l,\k))^+ \equiv_{forc} (P(\o )/\Fin )^+ \ast \pi$
and by Fact \ref{T209} the first iteration collapses $\fc ^V$ to $\fh ^V$.
Thus $\rp(\Col (\l,\k))^+ \Vdash |((\k^{<\l })^V)\check{\; }|\leq |(\fc ^V)\check{\;}|=|(\fh ^V )\check{\;}|$.
Clearly, statement (b) follows from (a).
\kdok
Under the assumption that $\cf (\l)=\o$ Theorem \ref{T233} gives sufficient conditions for the collapse of $\k ^{<\l}$ by $\rp (\Col(\l ,\k ))^+$.
In the sequel we show that $\cf (\l)\cf(\k ^{<\l})=\o$ implies that the cardinal $(\k ^{<\l})^+$ is collapsed too.
We use the following result of Baumgartner concerning almost disjoint families;
if $A(\k,\l,\mu,\nu)$ denotes the statement $\exists \CA \in [\k]^\mu \;(|\CA |=\l \land \forall A,B \in [\CA ]^2 \;|A\cap B|<\nu)$, we have
\begin{fac}\label{T353}
(Baumgartner \cite{Baum}) $A(\k,\k ^+,\cf (\k ),\cf (\k ))$, for each singular cardinal $\k$.
\end{fac}
\begin{te}\label{T356}
$(\ast)$ If $\cf (\l)\cf (\k^{<\l })=\o$, then

(a) $\rp (\Col(\l ,\k ))^+\Vdash |(((\k^{<\l }) ^+)^V)\check{\;}|\leq |(\fh ^V)\check{\;}|$;

(b) $\rp (\Col(\l ,\k ))^+\!\Vdash \!|(((\k^{<\l }) ^+)^V)\check{\,}|=\check{\o _1}$, if $\rp (\Col(\l ,\k ))^+\Vdash |(\fh ^V)\check{\,}|=\check{\o _1}$.
\end{te}
\dok
(a) Since $\cf (\l)=\o$, by Fact \ref{T351} we have $\Col (\l,\k)\cong \Col (\o ,\k^{<\l })$
and the statement follows from the following claim.
\begin{equation}\label{EQ355}
\k >\cf (\k)=\o \Rightarrow \rp (\Col(\o ,\k ))^+\Vdash |((\k ^+)^V)\check{\;}|\leq |(\fh ^V)\check{\;}|.
\end{equation}
We work with the preorder $\Q _{\o ,\k}$ and first show that $\rp (\Col(\o ,\k ))\Vdash |((\k ^+)^V)\check{\;}|\leq |(\fc ^V)\check{\;}|$.
In $V$,
the set $\CS :=\bigcup _{S\in [\o ]^\o} \{ S\}\times {}^{S}\o $ is the union of $\fc$-many sets of size $\fc$,
and, hence, $|\CS|=\fc$.
For $\p \in \Q _{\o ,\k}$ we have $S:=\supp(\p)\in [\o ]^\o$
and $\la \dom (\p  _i) :i\in \supp (\p)\ra \in {}^{S}\o$;
thus the function $f: \Q _{\o ,\k}\rightarrow \CS$ given by
$$
f(\p )=\Big\la \supp (\p), \la \dom (\p _i) :i\in \supp (\p)\ra\Big\ra, \;\;\mbox{ for }\p \in \Q _{\o ,\k},
$$
is well defined.
By Fact \ref{T353} $A(\k ,\k ^+,\o ,\o)$ is true;
so, there is a family of countable sets $\CA =\{ A^\xi : \xi < \k ^+\}\in [\k]^\o$
such that $|A^\xi \cap A^\zeta|<\o$, whenever $\xi <\zeta <\k ^+$.
For each $\xi <\k ^+$ we take an enumeration $A ^\xi=\{ \a ^\xi_i:i \in \o  \}$.
The sets
$$
\Delta _\xi :=\Big\{ \p\in Q _{\o ,\k} : \forall i\in \supp (\p ) \;\; \exists \eta _i \in {}^{<\o }\k\;\; \p _i =\eta _i{}^{\smallfrown}\a ^\xi_i\Big\}, \;\xi <\k ^+,
$$
are dense in $\Q _{\o ,\k}$.
Namely, if $\eta\in Q _{\o ,\k}$,
defining $\p _i=\eta _i{}^{\smallfrown}\a ^\xi_i$, for $i\in \supp (\eta )$, and $\p _i=s$, for $i\in \o \setminus \supp (\eta )$,
we have $\p =\la \p _i :i\in \o \ra\in \Delta_\xi$ and $\p\leq  \eta$.

Let $G$ be a $\Q _{\o ,\k}$-generic filter
and, in $V[G]$, let $\p^\xi \in G\cap \Delta _\xi$, for $\xi <\k ^+$.
Let $g:\k ^+ \rightarrow \CS$, where $g(\xi )=f(\p^\xi)=\la \supp (\p^\xi), \la \dom (\p ^\xi_i) :i\in \supp (\p^\xi)\ra\ra$, for $\xi <\k ^+$.
We show that $g$ is an injection, which will give $|\k^+|^{V[G]}\leq |\fc ^V|^{V[G]}$.

Otherwise, there would be $\xi <\zeta <\k ^+$ and $\la S,\la n_i +1 :i\in S\ra\ra \in \CS$
such that $F(\xi )=F(\zeta )=\la S,\la n_i +1 :i\in S\ra\ra $,
that is, $\supp (\p^\xi) =\supp (\p^\zeta)=S$ and $\dom (\p ^\xi_i ) = \dom (\p ^\zeta _i )=n_i +1 \in \o $, for $i\in S$.
So, for $i\in S$ we would have $\p ^\xi_i =\eta _i {}^{\smallfrown}\a ^\xi_i$ and $\p ^\zeta_i=\theta _i {}^{\smallfrown}\a ^\zeta_i $,
where $\eta _i =\p ^\xi_i \upharpoonright n_i $ and $\theta _i =\p ^\zeta_i \upharpoonright n_i$.
Since $\p^\xi , \p^\zeta \in G$ there is $\p\in G$ such that $\p\leq   \p^\xi , \p ^\zeta$
and by (\ref{EQ243}) there would be $j\in \o$ such that
$$
\forall i\in S\setminus j \;\;\p _i \supset \p^\xi_i , \p ^\zeta _i
$$
and, hence, $\a ^\xi_i =\p^\xi_i (n_i)=\p_i (n_i)=\p^\zeta_i (n_i)=\a ^\zeta_i$, for all $i\in S\setminus j$,
which is false because $|A^\xi \cap A^\zeta|<\o$. Thus, $|\k^+|^{V[G]}=|\fc ^V|^{V[G]}$.
Since $\rp(\Col (\o,\k))^+ \equiv_{forc} (P(\o )/\Fin )^+ \ast \pi$
and $(P(\o )/\Fin )^+$ collapses $\fc ^V$ to $\fh ^V$,
we have $\rp(\Col (\o,\k))^+ \Vdash |((\k^+)^V)\check{\; }|\leq |(\fc ^V)\check{\;}|=|(\fh ^V )\check{\;}|$
and (\ref{EQ355}) is proved.
Claim (b) follows from (a).
\hfill $\Box$
\subsection{Collapse of {\boldmath$\fh$} to {\boldmath$\o _1$}}
By Theorems \ref{T233} and  \ref{T356} the collapse of $\fh$ to $\o _1$ by $\rp (\Col (\l ,\k))$
provides the collapse of $\k ^{<\l}$ or $(\k ^{<\l})^+$ to $\o _1$, respectively.
That condition is trivially satisfied under $\fh =\o _1$
and here we detect more such conditions.
\begin{te}\label{T341}
$(\ast)$ If $\cf (\l )=\o$ and $\k ^{<\l} \geq \fc$, then $\rp(\Col (\l, \k ))^+ \Vdash |(\fc ^V)\check{\;}|=\check{\o _1}$.
\end{te}
\dok
By Fact \ref{T351} we have $\Col (\l, \k )\cong \Col (\o, \k ^{<\l})$;
so it is sufficient to prove that $\k \geq \fc$ implies $\rp(\Col (\o, \k ))^+ \Vdash |(\fc ^V)\check{\;}|=\o _1$
and then to replace $\k$ by $\k ^{<\l}$.

The poset $\P :={}^{<\o _1}\fc \times {}^{<\o}\k $ is a direct product of separative atomless posets;
so, it is separative and atomless.
By Fact \ref{T338}(e) we have ${}^{<\o }\k \hookrightarrow_c \P$ and, hence, $\P \Vdash |\k |=\o$.
In addition, $|P|= \fc ^\o \k =\k$
and by Theorem \ref{T321} (replacing $\l$ by $\o$) we have $\ro (\P )\cong \Col (\o, \k )$.
Now by Fact \ref{T338}(e) we have ${}^{<\o _1}\fc\hookrightarrow _c \P\hookrightarrow _d \ro (\P )^+ \cong \Col (\o , \k )^+$
so, by (a) and (b) of Fact \ref{T338}, ${}^{<\o _1}\fc\hookrightarrow _c  \Col (\o, \k )^+$,
and by Fact \ref{T338}(c) $\Col (\o _1, \fc)^+ =\ro ({}^{<\o _1}\fc)^+\hookrightarrow _c  \Col (\o, \k )^+$.
Finally, by Fact \ref{T300}(b) we have $\rp(\Col (\o _1, \fc))^+ \hookrightarrow _c  \rp(\Col (\o, \k ))^+$
and, by Theorem \ref{T340}(a), $\rp(\Col (\o _1, \fc))^+\equiv_{forc} \Col (\o _1,\fc)$;
thus, $\rp(\Col (\o, \k ))^+ \Vdash |(\fc ^V)\check{\;}|=\o _1$.
\hfill $\Box$

\begin{fac}\label{T220}
($0^\sharp$ does not exist)\\
If $\P$ is a separative $\o$-distributive forcing and $V_\P[G]$ a generic extension, then

(a) If $\fb =\fd $ and $\fb ^{V_\P[G]}=\o _1$, then $|\fb|^{V_\P[G]}=\o _1$;

(b) If $\cf ^{V_\P[G]}(\fh )=\o _1$, then $|\fh|^{V_\P[G]}=\o _1$.
\end{fac}
\dok
(a) In $V$, by $\fb =\fd $ there is a dominating family $\CD =\{ f_\a :\a <\fb\}\subset {}^\o\o$ well ordered by $\leq ^*$ (a scale)
and, since $\fb ^{V_\P[G]}=\o _1$, in $V_\P[G]$ there is an unbounded family $\{ g_\xi :\xi <\o _1\}$.
For $\xi <\o _1$ let $\a _\xi <\fb$, where $g_\xi \leq ^* f_{\a _\xi}$.
Assuming that the set $\{ a_\xi :\xi <\o _1\}$ is not cofinal in $\fb$
there would be $\a <\fb$ such that for each $\xi <\o _1$ we have $a_\xi \leq \a$ and, hence,   $g_\xi \leq ^* f_{\a _\xi}\leq ^* f_\a$,
which is false.
Thus the set $X:=\{ a_\xi :\xi <\o _1\}$ is cofinal in $\fb$.
Since $0^\sharp$ does not exist in $V_\P[G]$ too,
by Jensen's Covering Lemma applied in $V_\P[G]$ there is a set $C\in L^{V_\P[G]}=L \subset V$
such that $X\subset C$ and $|C|^{V_\P[G]}=\o _1$.
Clearly, the set $C\cap \fb$ is unbounded in $\fb$ and belongs to $V$
and, since $\fb$ is a regular cardinal (in $V$) we have $|C\cap \fb|^V =\fb$.
Hence we have $|\fb |^{V_\P[G]}=|C\cap \fb|^{V_\P[G]}\leq |C|^{V_\P[G]}=\o _1$
and, by the $\o$-distributivity of $\P$, $|\fb |^{V_\P[G]}\geq \o _1$,
which gives $|\fb|^{V_\P[G]}=\o _1$.

(b) Let $\cf ^{V_\P[G]}(\fh )=\o _1$ and, in $V_\P[G]$, let $X:=\{ a_\xi :\xi <\o _1\}$ be a cofinal subset of $\fh$.
Then, as above, there is a set $C\in L \subset V$
such that $X\subset C$ and $|C|^{V_\P[G]}=\o _1$.
Again, the set $C\cap \fh\in V$ is unbounded in $\fh$
and, since $\fh$ is a regular cardinal (in $V$) we have $|C\cap \fh|^V =\fh$.
Hence we have $|\fh |^{V_\P[G]}=|C\cap \fh|^{V_\P[G]}\leq |C|^{V_\P[G]}=\o _1$
and, by the $\o$-distributivity of $\P$, $|\fh |^{V_\P[G]}\geq \o _1$,
which gives $|\fh|^{V_\P[G]}=\o _1$.
\hfill $\Box$
\begin{te}\label{T228}
($\fb =\fd $ and $0^\sharp$ does not exist) \\
If $(\ast)$ holds and $\cf (\l)= \o $, then  $\rp(\Col (\l,\k))^+ \Vdash |(\fh ^V )\check{\;}|=\check{\o} _1$.
\end{te}
\dok
Note that by $(\ast)$ we have $\k ^{<\l } \geq \o _1 $.
Defining $\k _0=\k ^{<\l}$, by Fact \ref{T351} we have $\Col (\l,\k)\cong \Col (\o,\k _0)$;
by Proposition \ref{T301}(b) it is sufficient to prove that $\Q _{\o,\k _0}\Vdash |(\fh ^V )\check{\;}|=\check{\o}_1$.
First, it is evident that the sets
\begin{eqnarray*}
D_\xi  & := & \{ \p\in Q _{\o,\k _0} : \forall i\in \supp (\p )\;\; \xi \in \ran (\p _i)\}, \mbox{ for } \xi <\k_0,\\
E_S    & := & \{ \p\in Q _{\o,\k _0} : \supp (\p )\subset S \lor \supp (\p )\subset \o\setminus S\},  \mbox{ for } S\subset \o, \\
F_f    & := & \{ \p \in Q _{\o,\k _0}: \forall i\in \supp (\p ) \; |\p _i|\geq f(i) \},  \mbox{ for } f :\o\rightarrow \o ,
\end{eqnarray*}
are dense in $\Q _{\o,\k _0}$.
Let $V_{\Q _{\o,\k _0}}[G]$ be a generic extension;
suppose that $|\fh |^{V_{\Q _{\o,\k _0}}[G]}>\o _1$.
Then by Fact \ref{T220}(b) we have $\cf ^{V_{\Q _{\o,\k _0}}[G]}(\fh )>\o _1$.
Since $\fb \geq \fh$ we have $|\fb |^{V_{\Q _{\o,\k _0}}[G]}>\o _1$  and, by Fact \ref{T220}(a), $\fb ^{V_{\Q _{\o,\k _0}}[G]}>\o _1 $.

By Fact \ref{T300}(a) we have $\rp(\Col (\l,\k))^+ \equiv _{forc}(P(\o )/\Fin )^+ \ast \pi$
and (see \cite{Balc1}, p.\ 356) the generic filter in the first iteration is a selective ultrafilter $\CU $ in $P(\o )$
having a tower  of length $\fh $ as its base, say $\{ U_\a :\a <\fh\}$.
It is easy to show that $\{ \supp (\p): \p\in G \}\subset \CU $.

Let $\eta ^\xi \in G \cap D_\xi$, for $\xi <\o _1$.
Then for each $\xi <\o _1$ there is $\a _\xi \in \fh$ such that $U_{\a _\xi} \subset ^* \supp (\eta ^\xi) $.
Since $\cf ^{V_{\Q _{\o,\k _0}}[G]}(\fh )>\o _1$ we have $\b :=\sup\{ \a _\xi :\xi <\o _1\}<\fh$
and $U_{\b } \subset ^* U_{\a _\xi}\subset ^* \supp (\eta ^\xi)$, for $\xi <\o _1$.
Let $\eta ^{U_\b} \in G \cap E_{U_\b}$;
for $\xi <\o _1$ we have $\supp (\eta ^{U_\b} )\subset U_\b \subset ^* \supp (\eta ^\xi)$;
let $\f^\xi :\supp (\eta ^{U_\b}) \rightarrow \o$, where
\begin{equation}\label{EQ357}
\f^\xi (i) =       \left\{
                     \begin{array}{cl}
                         |\eta ^\xi_i| , & \mbox{ if } i\in \supp (\eta ^{U_\b}) \cap \supp (\eta ^\xi ),\\ [1mm]
                         0     , & \mbox{ if } i\in \supp (\eta ^{U_\b}) \setminus \supp (\eta ^\xi ).
                     \end{array}
                   \right.
\end{equation}
Since $\fb ^{V_{\Q _{\o,\k _0}}[G]}>\o _1 $ there is $\f :\supp (\eta ^{U_\b}) \rightarrow \o$ such that $\f^\xi \leq ^* \f$, for all $\xi <\o _1$.
Let $\f \subset f :\o\rightarrow \o$. Then for each $\xi <\o _1$ there is $k_\xi \in \o$ such that
\begin{equation}\label{EQ255}
\forall i\in \supp (\eta ^{U_\b}) \cap \supp (\eta ^\xi )\setminus k_\xi \;\; |\eta ^\xi_i|\leq f(i).
\end{equation}
Let $\eta ^f\in G \cap F_f$.
Since $\eta ^{U_\b},\eta ^f \in G$ there is $\p \in G$ such that $\p \leq _{Q _{\o,\k _0}} \eta ^{U_\b},\eta ^f$
and by (\ref{EQ243}) there is $j^*\in \o$ such that
\begin{equation}\label{EQ244}
\forall i\in \supp (\p)\setminus j^* \;\;\Big(i\in \supp (\eta ^{U_\b})\cap \supp (\eta ^f) \land \eta^{U_\b}_i,\eta ^f_i \subset \p _i\Big ).
\end{equation}
We will show that $\o _1 \subset \bigcup _{i\in \supp (\p)}\ran (\p_i)=:R_\p$,
which is impossible since the set $R_\p$ is countable and $\o _1$ is not collapsed.
Let $\xi <\o _1$. Since $\eta ^\xi \in G$ there is $\p ^\xi \in G$ such that $\p ^\xi \leq \eta ^\xi ,\p$.
By (\ref{EQ243}) there is $j_\xi\in \o$ such that
\begin{equation}\label{EQ254}
\forall i\in \supp (\p ^\xi)\setminus j_\xi \;\;\Big(i\in \supp (\eta ^\xi )\cap \supp (\p ) \land \eta^\xi _i,\p _i \subset \p ^\xi _i\Big ).
\end{equation}
Let $j:=\max\{k_\xi , j^*, j_\xi \}$ and let $i\in \supp (\p ^\xi)\setminus j$.
By (\ref{EQ254}) we have $i\in \supp (\eta ^\xi )$ and $i\in \supp (\p )\setminus j^*$,
which  by (\ref{EQ244}) gives $i\in \supp (\eta ^{U_\b})\cap \supp (\eta ^f)$; thus,
\begin{equation}\label{EQ256}
i\in (\supp (\p ^\xi) \cap\supp (\eta ^{U_\b})\cap \supp (\eta ^f)\cap \supp (\eta ^\xi )\cap \supp (\p ))\setminus j.
\end{equation}
By (\ref{EQ256}), (\ref{EQ255}), and since $\eta ^f \in F_f$ we have $|\eta ^\xi_i|\leq f(i)\leq |\eta ^f_i|$.
By (\ref{EQ244}) and (\ref{EQ254}) we have $\eta ^f_i \subset \p _i\subset \p ^\xi _i$ and $\eta^\xi _i \subset \p ^\xi _i$;
so,  $\eta^\xi _i$ and $\eta ^f_i$ are comparable,
which, together with $|\eta ^\xi_i|\leq |\eta ^f_i|$, gives $\eta^\xi _i\subset \eta ^f_i\subset \p _i$.
Since $\eta ^\xi \in D_\xi$ by (\ref{EQ256}) we have $\xi \in \ran (\eta ^\xi)\subset \ran (\p _i)\subset R_\p$.
\hfill $\Box$
\subsection{Boolean completions}
\begin{te}\label{T360}
($\ast$) Let $\rp (\Col (\l ,\k))^+\Vdash |(((\k ^{<\l })^\o)^V)\check{\;}|=\check{\o _1}$ and
\begin{equation}\label{EQ335}
\cf (\l)\cf(\k ^{<\l})=\o \;\mbox{ or }\; \l =\o .
\end{equation}
Then for each $n\in \N$ we have
\begin{equation}\label{EQ337}
\ro (\rp ^n (\Col (\l ,\k)))\cong \Col (\o _1 ,(\k ^{<\l})^\o).
\end{equation}
\end{te}
\dok
Let $\f (\l,\k)$ denote the conjunction of all assumptions on $\l$ and $\k$. By induction we prove that for each $n\in \N$ we have
\begin{equation}\label{EQ336}
\forall \l,\k\in \Card \;\;\Big(\f (\l ,\k)\Rightarrow \ro (\rp ^n (\Col (\l ,\k)))\cong \Col (\o _1 ,(\k ^{<\l})^\o)\Big).
\end{equation}
First let $n=1$ and let $\l$ and $\k$ be cardinals satisfying  $\f (\l ,\k)$.
Since (\ref{EQ335}) implies that $\cf (\l)=\o$,
by Fact \ref{T351} and Proposition \ref{T301}(b) we have $\ro (\rp (\Col (\l ,\k)))\cong\ro (\rp (\Col (\o,\k ^{<\l})))\cong \ro (\sq (\Q _{\o,\k ^{<\l}}))$;
so, by Fact \ref{T2226}(a), $\rp (\Col (\l ,\k))^+ \equiv_{forc}\Q _{\o,\k ^{<\l}}$
and, hence, $\Q _{\o, \k ^{<\l }}\Vdash |(((\k ^{<\l })^\o)^V)\check{\;}|=\check{\o _1}$.

For a proof that $\ro (\sq (\Q _{\o,\k ^{<\l}}))\cong \Col (\o _1 ,(\k ^{<\l})^\o)$ we use Corollary \ref{T365}
in which we replace $\mu$ by $(\k ^{<\l})^\o$ and $\P$ by $\Q _{\o,\k ^{<\l}}$.
So, by Proposition \ref{T301}(a)
and since $(\k ^{<\l})^\o =((\k ^{<\l})^\o)^\o \geq \o _1$ and $\Q _{\o, \k ^{<\l }}\Vdash |(((\k ^{<\l })^\o)^V)\check{\;}|=\check{\o _1}$,
it remains to be proved that $\Q _{\o, \k ^{<\l }} =(\k ^{<\l})^\o$
First, let $\cf(\k ^{<\l})=\o$.
From $\k\l > \o$ it follows that $\k ^{<\l}\geq \o _1$;
so we have $(\k ^{<\l})^\o=(\k ^{<\l})^{\cf (\k ^{<\l})}> \k ^{<\l}\geq \o _1$.
Since $\Q _{\o,\k ^{<\l}}$ collapses $(\k ^{<\l})^\o>\o_1$ to $\o _1$, we have $|\Q _{\o,\k ^{<\l}}|\geq(\k ^{<\l})^\o$;
on the other hand, by Proposition \ref{T301}(a), $|\Q _{\o,\k ^{<\l}}|\leq ((\k ^{<\l})^{<\o})^\o=(\k ^{<\l})^\o$,
thus $|\Q _{\o,\k ^{<\l}}|=(\k ^{<\l})^\o $,
and we are done.
Second, if $\l =\o$,
then $\k \geq\o _1$ and we have to prove that  $|\Q _{\o,\k }|=\k ^\o$.
By Proposition \ref{T301}(a) we have $\k \leq|\Q _{\o,\k }|\leq\k ^\o$; so, if  $\k ^\o =\o _1$, we are done.
If $\k ^\o >\o _1$, then, since $\Q _{\o,\k }$ collapses $\k ^\o$ to $\o _1$, we have $|\Q _{\o,\k }|\geq\k ^\o$ and (\ref{EQ336}) is true for $n=1$.

Assuming that (\ref{EQ336}) is true (for $n$) we prove that it is true for $n+1$.
Let $\l$ and $\k$ be cardinals satisfying  $\f (\l ,\k)$.
Then by our assumption we have (\ref{EQ337})
and first prove that
\begin{equation}\label{EQ338}
\ro( \rp (\Col (\o _1, (\k ^{<\l})^\o)))\cong  \Col (\o _1, (\k ^{<\l})^\o).
\end{equation}
We regard Theorem \ref{T359} and first replace $\l$ by $\o _1$ and take $n=1$.
The second condition in (\ref{EQ331}) is satisfied; so for each $\mu \geq 2$
from $\;\rp (\Col (\o _1 ,\mu))^+\Vdash |((\mu ^{\o })^V)\check{\;}|=\check{\o _1}$
it follows that $\ro (\rp (\Col (\o _1 ,\mu)))\cong \Col (\o _1 ,\mu ^{\o})$.
Further, replacing $\mu$ by $(\k ^{<\l})^\o$, since $((\k ^{<\l})^\o) ^{\o}=(\k ^{<\l})^\o$, we conclude that
$\;\rp (\Col (\o _1 ,(\k ^{<\l})^\o))^+\Vdash |(((\k ^{<\l}) ^{\o })^V)\check{\;}|=\check{\o _1}$
implies $\ro (\rp (\Col (\o _1 ,(\k ^{<\l})^\o)))\cong \Col (\o _1 ,(\k ^{<\l}) ^{\o})$, that is,  (\ref{EQ338}).
Thus, it remains to prove that $\;\rp (\Col (\o _1 ,(\k ^{<\l})^\o))^+\Vdash |(((\k ^{<\l}) ^{\o })^V)\check{\;}|=\check{\o _1}$.
Replacing $\l$ by $\o _1$ in (a) of Theorem \ref{T352}
we conclude that for each $\mu \geq 2$ we have
$\rp(\Col (\o _1,\mu))^+\Vdash |((\mu ^{\o })^V)\check{\;}|=\check{\o _1}$.
Thus for $\mu =(\k ^{<\l})^\o$
we have $\rp(\Col (\o _1,(\k ^{<\l})^\o))^+\Vdash |(((\k ^{<\l}) ^{\o })^V)\check{\;}|=\check{\o _1}$ and (\ref{EQ338}) is true.
Now we have
\begin{eqnarray*}
\ro(\rp ^{n+1}(\Col (\l,\k))) & =     & \ro( \rp (\rp ^n(\Col (\l,\k))))        \quad \hfill (\mbox{by definition of }\rp^n)\\
                              & \cong & \ro( \rp (\ro (\rp ^n(\Col (\l,\k)))))  \quad \hfill (\mbox{by Fact \ref{T334}(a)})\\
                              & \cong & \ro( \rp (\Col (\o _1, (\k ^{<\l})^\o)))     \quad \hfill (\mbox{by (\ref{EQ337})})\\
                              & \cong & \Col (\o _1, (\k ^{<\l})^\o)                 \quad \hfill (\mbox{by (\ref{EQ338})}) \quad \Box.
\end{eqnarray*}
In the sequel we show that in many models of ZFC all Boolean algebras of the form $\ro (\rp ^n(\Col (\l,\k)))$, where $\cf (\l )=\o$, are collapsing algebras.
We recall that the Singular Cardinal Hypothesis (SCH) is
the statement that for each singular cardinal $\k$ we have
\footnote{We recall that the theories ZFC + $\neg$ SCH  and ZFC + ``There is a measurable cardinal $\k$ of Mitchell order $\k ^{++}$" are equiconsistent  \cite{Git}
and if $0^\sharp$ does not exist, then SCH holds (see \cite{Jech}, p.\ 357). In addition, SCH follows from GCH and PFA.}
\begin{equation}\label{EQ350}
\k ^{\cf (\k)}=\max \{ 2^{\cf (\k)},\k ^+\}.
\end{equation}
\begin{fac}\label{T363}
$(\ast)$ Under SCH we have
\begin{equation}\label{EQ348}
 (\k^{<\l}) ^\o = \left\{
                     \begin{array}{ll}
                          \fc  ,     & \mbox{ if } \k^{<\l} \leq \fc ;\\[1mm]
                          \k^{<\l}   ,     & \mbox{ if } \k^{<\l}    > \fc \land \cf (\k^{<\l})>\o;\\[1mm]
                          (\k^{<\l}) ^+,     & \mbox{ if } \k ^{<\l}   > \fc \land \cf (\k^{<\l})=\o.
                     \end{array}
         \right.
\end{equation}
\end{fac}
\dok
We recall that, in ZFC, for a cardinal $\k>\fc$ we have (see \cite{Jech}, p.\ 49)
\begin{equation}\label{EQ349}
\k ^\o = \left\{
                     \begin{array}{ll}
                          \mu ^\o,        & \mbox{ if }  \exists \mu <\k \; \mu ^\o \geq \k ;\\[1mm]
                          \k    ,         & \mbox{ if }  \forall \mu <\k \; (\mu ^\o < \k) \land \cf (\k)>\o;\\[1mm]
                          \k ^{\cf (\k)}, & \mbox{ if }  \forall \mu <\k \; (\mu ^\o < \k) \land \cf (\k)=\o.
                     \end{array}
         \right.
\end{equation}
For a proof of (\ref{EQ348}) it is sufficient to show that under SCH for each cardinal $\k \geq 2$ we have
\begin{equation}\label{EQ346}
\k ^\o = \left\{
                     \begin{array}{ll}
                          \fc  ,     & \mbox{ if } \k \leq \fc ;\\[1mm]
                          \k   ,     & \mbox{ if } \k    > \fc \land \cf (\k)>\o;\\[1mm]
                          \k ^+,     & \mbox{ if } \k    > \fc \land \cf (\k)=\o.
                     \end{array}
         \right.
\end{equation}
Clearly, $2\leq \k \leq \fc$ gives $\k ^\o =\fc$.
If $\k > \fc$ and $\cf (\k)=\o$,
then by (\ref{EQ350}) $\k ^\o =\k ^{\cf (\k )}\in \{ \fc , \k ^+\}$
and, since $\k ^\o \geq \k >\fc$, we obtain $\k ^\o =\k ^+$.

Finally let $\k > \fc$ and $\cf (\k)>\o$.
Assuming towards a contradiction that $\k ^\o >\k$
by (\ref{EQ349}) there would be $\mu <\k$ such that $\mu ^\o \geq \k$;
let $\mu$ be the minimal cardinal with these properties.
Since $\mu \leq \fc$ would imply $\k \leq\mu ^\o =\fc$ we have $\mu >\fc$ and, by (\ref{EQ349}),
\begin{equation}\label{EQ347}
\mu ^\o = \left\{
                     \begin{array}{ll}
                          \nu ^\o,          & \mbox{ if } \exists \nu <\mu \; \nu ^\o \geq \mu ;\\[1mm]
                          \mu    ,          & \mbox{ if } \forall \nu <\mu \; (\nu ^\o < \mu) \land \cf (\mu)>\o;\\[1mm]
                          \mu ^{\cf (\mu)}, & \mbox{ if } \forall \nu <\mu \; (\nu ^\o < \mu) \land \cf (\mu)=\o.
                     \end{array}
         \right.
\end{equation}
Assuming that there is $\nu <\mu\leq \nu ^\o $
we would have $\nu ^\o\leq \mu^\o\leq \nu ^\o $
and, hence, $\nu ^\o = \mu^\o \geq \k$, which is false by the minimality of $\mu$.
Thus by (\ref{EQ347}) and since $\mu ^\o >\mu$ we have $\cf (\mu)=\o$.
So, by (\ref{EQ350}), $\mu^\o =\mu ^{\cf (\mu)}\in \{ \fc , \mu ^+\}$ and, since $\mu >\fc$, we have $\mu^\o =\mu ^+$.
But this is false because
$\mu <\k \leq \mu ^\o$
and, hence, $\mu^\o =\k^\o >\k >\mu$
and we have a contradiction.
\hfill $\Box$
\begin{te}\label{T335}
(a) If $(\ast)$ holds, $\cf (\l )=\o$ and $\rp (\Col (\l,\k))^+ \Vdash |(\fh ^V )\check{\;}|=\check{\o_1}$, then for each $n\in \N$
\begin{equation}\label{EQ340}
\ro (\rp ^n(\Col (\l,\k)))\cong  \left\{
                                        \begin{array}{ll}
                                           \Col (\o _1, \k ^{<\l }),    & \mbox{ if } (\k ^{<\l})^\o =\k ^{<\l};                                                  \\ [1mm]
                                           \Col(\o _1,(\k^{<\l }) ^+ ), & \mbox{ if } \cf (\k^{<\l })=\o \;\mbox{ and }\; (\k^{<\l }) ^\o=(\k^{<\l }) ^+; \\ [1mm]
                                           \Col (\o _1, \fc ),          & \mbox{ if } \l =\o             \;\mbox{ and }\; \k\in [\o _1, \fc].
                                        \end{array}
                                 \right.
\end{equation}

(b) Under SCH formula (\ref{EQ340}) determines $\ro (\rp ^n(\Col (\l,\k)))$
for each pair of cardinals $\l,\k$ satisfying the assumptions from (a).

(c) Statement (a) is true if instead of $\rp (\Col (\l,\k))^+ \Vdash |(\fh ^V )\check{\;}|=\check{\o_1}$ 
we have
\begin{equation}\label{EQ358}
\k ^{<\l } \geq \fc \;\mbox{ or }\; \fh =\o _1 \;\mbox{ or }\; (\fb =\fd \;\mbox{ and }\;0^\sharp \mbox{ does not exist}).
\end{equation}
\end{te}
\dok
(a) Assuming that $\cf (\l )=\o$ and $(\k ^{<\l})^\o =\k ^{<\l}$
we show that the assumptions of Theorem \ref{T359} are satisfied.
By Theorem \ref{T233}(b) we have $\rp(\Col (\l,\k))^+ \Vdash |((\k ^{<\l })^V)\check{\; }|= \check{\o} _1$
and we show that (\ref{EQ331}) is true.
If $\l =\o$, then by $(\ast)$ we have $\k >\o$ and, since $(\k ^{<\l})^\o =\k ^{<\l}$ and $\k^{<\o}=\k$, we have $\k ^\o =\k $ and (\ref{EQ331}) holds.
If $\l >\o$, then assuming that $\cf (\k ^{<\l})=\o$
we would have  $(\k ^{<\l})^\o =(\k ^{<\l})^{\cf (\k ^{<\l})}>\k ^{<\l}$, which contradicts our assumption.
Thus  $\cf (\k ^{<\l})>\o$ and (\ref{EQ331}) holds again.

Let $\cf (\l )\cf (\k^{<\l })=\o$ and $(\k^{<\l }) ^\o=(\k^{<\l }) ^+$.
By Theorem \ref{T356}(b) and since $(\k^{<\l }) ^\o=(\k^{<\l }) ^+$ we have $\rp (\Col(\l ,\k ))^+\Vdash |(((\k^{<\l }) ^\o)^V)\check{\;}|=\check{\o _1}$.
Since  $\cf (\l)\cf (\k^{<\l })=\o$, the statement follows from Theorem \ref{T360}.

Finally, let $\l =\o$ and $\k\in [\o _1, \fc]$.
We use Theorem \ref{T360} for $\l =\o$. Since $\k\in [\o _1, \fc]$ we have $(\k ^{<\l })^\o=\k ^\o =\fc $;
thus, $\rp (\Col (\o ,\k))^+\Vdash |(\fc ^V)\check{\;}|=\check{\o _1}$ implies that $\ro (\rp ^n(\Col (\o,\k)))\cong  \Col (\o _1, \fc )$ holds
for each $n\in \N$.
But by Fact \ref{T300}(a) we have $\rp (\Col (\o,\k))^+ \equiv_{forc}\CP _\o \ast \pi$
and, by Fact \ref{T209}, $\rp (\Col (\o,\k))^+ \Vdash |(\fc ^V )\check{\;}|=|(\fh ^V )\check{\;}|$,
which together with $\rp (\Col (\o,\k))^+ \Vdash |(\fh ^V )\check{\;}|=\check{\o_1}$ gives
$\rp (\Col (\o,\k))^+ \Vdash |(\fc ^V )\check{\;}|=\check{\o_1}$ and we are done.

(b) Assuming that the first two conditions in (\ref{EQ340}) are not satisfied we show that the third holds.
First, if $(\k ^{<\l})^\o >\k ^{<\l}$ and $\cf (\k^{<\l })>\o$,
then by (\ref{EQ348}) from Fact \ref{T363} we have $(\k ^{<\l})^\o =\fc$.
Thus $\k ^{<\l}<\fc$ and, hence, $\l =\o$ and $\k <\fc$ and the third condition in (\ref{EQ340}) is satisfied.
Second, if  $(\k ^{<\l})^\o >\k ^{<\l}$ and $(\k^{<\l }) ^\o\neq(\k^{<\l }) ^+$,
then, by (\ref{EQ348}), $(\k ^{<\l})^\o =\fc$.
So $\k ^{<\l}<\fc$ again and we continue as above.

(c) follows from Theorems \ref{T341} and \ref{T228}.
\hfill $\Box$
\section{{\boldmath$\l $}'s of cofinality {\boldmath$\leq \fc $}}\label{S7}
\begin{te}\label{T364}
{\rm (SCH + $\fh =\o _1$)} If $(\ast)$ holds and $\cf (\l )\leq \fc$, then for each $n\in \N$ we have
\begin{equation}\label{EQ352}
\ro (\rp ^n(\Col (\l,\k)))\cong \Col (\o _1, (\k ^{<\l })^\o );
\end{equation}
more precisely,
\begin{equation}\label{EQ351}
\ro (\rp ^n(\Col (\l,\k)))\cong   \left\{
                                         \begin{array}{ll}
                                              \Col (\o _1,\fc )          ,     & \mbox{ if } \k^{<\l} \leq \fc ;\\[1mm]
                                              \Col (\o _1,\k^{<\l})      ,     & \mbox{ if } \k^{<\l}    > \fc \land \cf (\k^{<\l})>\o;\\[1mm]
                                              \Col (\o _1,(\k^{<\l}) ^+ ),     & \mbox{ if } \k ^{<\l}   > \fc \land \cf (\k^{<\l})=\o.
                                         \end{array}
                                  \right.
\end{equation}
If instead of {\rm (SCH + $\fh =\o _1$)} we have $( 0^\sharp \mbox{ does not exist}+ \fb =\fd)$, then the previous statement holds when $\cf (\l )= \o$.
\end{te}
\dok
If $\cf (\l) \in [\o _1,\fc]$,
then, by Theorem \ref{T340}(a) we have $\ro (\rp ^n (\Col (\l ,\k)))\cong \Col (\o _1 ,\k ^{<\l})$, for all $n\in \N$.
Since $\o <\cf (\l)$ by Fact \ref{T327}(a) we have $(\k ^{<\l })^\o \leq (\k ^{<\l })^{<\cf (\l)}=\k ^{<\l }\leq (\k ^{<\l })^\o$;
thus, $(\k ^{<\l })^\o=\k ^{<\l}$ and we are done.

If $\cf (\l )=\o$, then since $\fh =\o _1$ by Theorem \ref{T335}(a) we have (\ref{EQ340})
and by Theorem \ref{T335}(b) at least one condition from (\ref{EQ340}) is satisfied by $\l $ and $\k$.
First, if $(\k ^{<\l})^\o =\k ^{<\l}$,
then by (\ref{EQ340}) we have $\ro (\rp ^n(\Col (\l,\k)))\cong\Col (\o _1, \k ^{<\l })=\Col (\o _1, (\k ^{<\l })^\o )$
and is (\ref{EQ352}) true.

Otherwise we have $(\k ^{<\l})^\o >\k ^{<\l}$
and regarding (\ref{EQ348}) there are two cases.

{\it Case 1}: $\k ^{<\l}   > \fc$,  $\cf (\k^{<\l})=\o$ and $(\k^{<\l }) ^\o=(\k^{<\l }) ^+$.
Then by (\ref{EQ340}) we have $\ro (\rp ^n(\Col (\l,\k)))\cong\Col (\o _1, (\k ^{<\l })^+)=\Col (\o _1, (\k ^{<\l })^\o )$
and is (\ref{EQ352}) true.

{\it Case 2}: $\k ^{<\l}   \leq \fc$ and $(\k^{<\l})^\o=\fc$.
Then, since $(\k ^{<\l})^\o >\k ^{<\l}$,
we have $\k ^{<\l}   < \fc$,
which by $(\ast)$ implies that $\l =\o$ and $\k \in [\o _1, \fc )$.
So, by (\ref{EQ340}) we have $\ro (\rp ^n(\Col (\l,\k)))\cong\Col (\o _1, \fc )=\Col (\o _1, (\k ^{<\l })^\o )$
and is (\ref{EQ352}) true again.

Statement (\ref{EQ351}) follows from (\ref{EQ352}) and Fact \ref{T363}.

If $0^\sharp \mbox{ does not exist}$ and $\fb =\fd$, then we have SCH and if $\cf (\l)=\o$ by Theorem \ref{T335}(a) we have (\ref{EQ340}).
Now we prove (\ref{EQ352}) as above.
\hfill $\Box$
\begin{ex}\label{EX302}\rm
Under GCH, if $(\ast)$ holds and $\cf (\l )\leq \o _1$, then by (\ref{EQ352}) from Theorem \ref{T364} and Fact \ref{T327}(b) we have

 - $\ro (\rp ^n(\Col (\l,\k)))\cong \Col (\o _1,\mu)$, for all $n\in \N$, where $\mu =\k^{<\l}$, if $\cf (\k^{<\l})>\o$, and $\mu =(\k^{<\l}) ^+ $, if $\cf (\k^{<\l})=\o$;

 - $\ro(\rp ^n (\Col (\l ,\k)))\cong \Col (\o _1,\mu)$, if  $\l \leq \o _1\leq \k $ and $n\in \N$, where $\mu =\k$, if $\cf (\k)>\o$, and $\mu =\k ^+$, if $\cf (\k)=\o$.

\noindent
More generally, by Theorem \ref{T340}(a) we have $\ro(\rp ^n (\Col (\o _1,\k)))\cong \Col (\o _1, \k ^{\o})$, for all $n\in \N$.
So, by Theorem \ref{T335}, if $\fh =\o _1$ or ($\fb =\fd $ and $0^\sharp$ does not exist) and if $\k =\k ^\o $,
then $\ro (\rp ^n(\Col (\o,\k)))\cong \ro(\rp ^n (\Col (\o _1,\k)))\cong \Col (\o _1, \k )$, for all $n\in \N$.
\end{ex}
\begin{ex}\label{EX303}\rm
If $0^\sharp$ does not exist and MA + $\fc =\o _2$ holds, then SCH and $\fp=\ft=\fh=\fb =\fd=\fc=\o _2$ hold too.
By Theorems \ref{T340}(a) and \ref{T335} for $\l \in \{ \o ,\o _1\}$ and $\k \in \{ \o _1 ,\o _2\}$ we have $\ro(\rp ^n (\Col (\l ,\k)))\cong \Col (\o _1,\o _2)$.
But $\h (\rp (\Col (\o ,\o)))\geq\ft =\o _2$ (Balcar and Hru\v{s}\'{a}k, \cite{BH}) and $\h (\rp (\Col (\o _2 ,\o _2)))=\o _2$ (Dow, \cite{Dow1});
thus the reduced powers $\rp (\Col (\o ,\o))$ and $\rp (\Col (\o _2 ,\o _2))$ do not collapse cardinals.
\end{ex}
\begin{ex}\label{EX304}\rm
If $V$ is a model of GCH and $\P =\Fn (\o_2 ,2)$ is the Cohen forcing and $\mu$ and $\theta$ are infinite cardinals, then
$\fc ^{V[G]}=\o_2$ and $(\theta ^\mu)^{V[G]}=(\max(\o _2 ,\theta)^\mu)^V $ (see \cite{Kun}, p.\ 246).
In $V[G]$ we have $\fh =\o _1$ (see \cite{Blas})
and if $\l \in \{ \o ,\o _1,\o _2\}$ and $\k \in \{ \o _1 ,\o _2\}$, then $\k ^{<\l }\leq \o _2 ^{\o _1}=\o _2=\fc$
and by Theorem \ref{T364} and Fact \ref{T367} $\ro (\rp ^n(\Col (\l,\k)))\cong \Col (\o _1, \fc)\cong \ro (P(\o)/\Fin)$.
\end{ex}
\section{Some extensions}\label{S8}
\paragraph{Incarnations of collapsing algebras}
Collapsing algebras have several incarnations;
for example, some basic examples of Boolean algebras are, in fact, collapsing algebras:
$\RO (2^\o)\cong\Col (\o ,2)$ and $\ro (P(\o )/\Fin)\cong\Col (\o _1,2)$ under CH.
In addition, Boolean completions of reduced powers of some algebras are isomorphic to regular open algebras of some spaces of ultrafilters.
These standard facts are generalized in the forthcoming Fact \ref{T362},
providing alternative interpretations of the results obtained in the previous sections;
first we recall some definitions.

For a cardinal $\k\geq \o$, by $\CP _\k $ we denote the partial order $(P(\k )/[\k ]^{<\k })^+$
and by $U_\k$ the space of uniform ultrafilters in $P(\k )$.

If $\CX=( X,\CO )$ is a topological space,
by $\Clop (\CX)$ we denote the Boolean algebra of clopen subsets of $\CX$ and
by $\RO (\CX)$ the complete Boolean algebra of regular open subsets of $\CX$.

If $\l \geq \o$ is a regular cardinal,
a space $\CX$ is called a {\it $P_\l$-space}, or {\it $\l$-open},
iff the intersection of $<\l$-many open sets is open (possibly empty).
$\CX _\l$ denotes the space $( X,\CO_\l )$, where $\CO_\l$ is the smallest $\l$-open topology on the set $X$ containing $\CO$.
$(\CX^\l)_\l$ denotes the box-product of $\l$-many copies of $\CX$ with boxes of size $<\l$;
if, in particular, $\CX$ is a cardinal $\k$ with the discrete topology, we write $(\k^\l)_\l$.
\begin{fac}\label{T362}
Let $\l\geq \o$ and $\k \geq 2$ be cardinals. Then

(a) $\RO ((\k^\l)_\l)\cong \Col(\l ,\k)$; in particular,  $\RO ((2^\l)_\l)\cong\Col(\l ,2)$, if $\l$ is a regular cardinal;

(b) $ \ro (\CP _\k)\cong\RO (U_\k)$, if $\k \geq \o$;

(c) If $2^\k =\k ^+$ and $2^{\cf(\k )} =\cf (\k )^+$, then
\begin{equation}\label{EQ257}
\ro (\CP _\k )\cong \left\{
                     \begin{array}{cl}
                         \Col (\o \,,\k ^+) , & \mbox{ if }\cf (\k )>\o ,\\ [1mm]
                         \Col (\o _1,\k^+), & \mbox{ if }\cf (\k )=\o ;
                     \end{array}
                    \right.
\end{equation}

(d) $\RO ((U_\k )_{2^\k})\cong \RO((2^{(2^\k)})_{2^\k})\cong\Col(2^\k ,2)$, if $2^\k $ is regular and each base for each ultrafilter $\CU \in U_\k$ has size $2^\k$;

(e) $\RO ((U_\k )_{\k ^+})\cong \RO ((2^{(\k ^+)})_{\k^+})\cong\Col(\k ^+ ,2)$, if $2^\k =\k ^+$.
\end{fac}
\dok
(a) For $\f \in {}^{<\l}\k$ let $B_\f :=\{ x\in \k^\l: \f \subset x\}$.
Clearly, the set $\CB :=\{ B_\f :\f \in {}^{<\l}\k \}$ is a clopen base for the topology on $(\k^\l)_\l$,
and ${}^{<\l}\k $ densely embeds in $\Clop ((\k^\l)_\l )$ which densely embeds in $\RO ((\k^\l)_\l)$.
Since $\Col(\l ,\k) =\ro ({}^{<\l}\k)$ we have $\Col(\l ,\k)\cong \RO ((\k^\l)_\l)$.

Concerning (b), by \cite{Comf} (Lemma 7.12, p.\ 151) we have $\CP_\k\cong\Clop(U_\k)$ and, hence, $ \ro (\CP _\k)\cong\RO (U_\k)$.
For (c) see \cite{Balc1}, p.\ 380.
For a proof of the first part of (d) and (e) see \cite{Comf}, p.\ 401; the second part follows from (a).
\kdok
\noindent
So, for example, we obtain the following consequence of Theorems \ref{T335} and \ref{T340}.
\begin{te}\label{T344}
Let $\k\geq \o$ be a cardinal,  $2^\k =\k ^+$ and $2^{\cf (\k )}=\cf (\k )^+$.

(a) ($\h =\o _1$ or ($\fb=\fd$ and $0 ^\sharp$ does not exist)) If  $\cf (\k)>\o$, then $\Col (\o ,\k ^+)\cong  \RO ((\k ^+)^\o)\cong\ro (\CP _\k)\cong \RO (U_\k)$
and if $\B$ is any of these algebras or $\B =\CP _\k$,
then for each $n\in \N$ we have
\begin{equation}\label{EQ323}
\ro(\rp ^n (\B))\cong \Col (\o _1, \k ^+ );
\end{equation}

(b)  If $\cf (\k )=\o$, then $\Col (\o _1 ,\k ^+)\cong \RO (((\k ^+)^{\o _1})_{\o _1})\cong \ro (\CP _\k)\cong \RO (U_\k)$
and if $\B$ is any of these algebras or $\B =\CP _\k$,
then for each $n\in \N$ we have
\begin{equation}\label{EQ324}
\ro(\rp ^n (\B))\cong \ro(\B)\cong \Col (\o _1, \k ^+ ).
\end{equation}
\end{te}
\dok
(a) First we use Fact \ref{T362}. By (a) we have $\Col(\o ,\k ^+)\cong \RO (((\k^+)^\o)_\o)=\RO ((\k^+)^\o)$;
by (c), $\ro (\CP _\k)\cong\Col (\o ,\k ^+)$;
and by (b) $\ro (\CP _\k)\cong \RO (U_\k)$.
Now, if $\B\neq\CP _\k$, then, since $\B \cong \Col (\o ,\k^+)$,
we have $\ro(\rp ^n (\B))\cong \ro(\rp ^n (\Col (\o ,\k^+)))$
and we use Theorem \ref{T335}, where we replace $\l$ by $\o$ and $\k$ by $\k ^+$.
Since $((\k ^+)^{<\o})^\o =2^{\k\o}=2^\k=\k ^+=(\k ^+)^{<\o}$,
the first condition in (\ref{EQ340}) is satisfied;
so we obtain $\ro(\rp ^n (\Col (\o ,\k^+)))\cong \Col (\o _1, \k ^+ )$, for each $n\in \N$, and (\ref{EQ323}) is true.
If $\B=\CP _\k$,
then by Fact \ref{T334}(b) we have $ \ro(\rp ^n (\Col (\o ,\k ^+)))\cong  \ro(\rp ^n (\ro(\CP _\k)))\cong  \ro(\rp ^n (\CP _\k))$
and (\ref{EQ323}) is true again.

(b) Again we use Fact \ref{T362}.
By (a) we have $\Col(\o_1 ,\k ^+)\cong \RO (((\k^+)^{\o_1})_{\o_1})$;
by (c), $\ro (\CP _\k)\cong\Col (\o _1 ,\k ^+)$;
and by (b) $\ro (\CP _\k)\cong \RO (U_\k)$.
Now,  if $\B\neq\CP _\k$,
then, since $\B\cong \Col (\o _1 ,\k ^+)$,
we have $\ro(\rp ^n (\B))\cong \ro(\rp ^n (\Col (\o _1 ,\k^+)))$
and we use Theorem  \ref{T340}(b), in which we replace $\l$ by $\o_1$ and $\k$ by $\k ^+$.
(Note that by the assumptions we have $2^\o =\o _1$ and, hence, $\h =\o _1$.)
Since $(\k ^+)^{<\o _1} =2^{\k\o}=2^\k=\k ^+$
we have $\ro(\rp ^n (\Col (\o _1 ,\k^+)))\cong\Col (\o _1, \k ^+ )$, for all $n\in \N$, and (\ref{EQ324}) is true.
If $\B=\CP _\k$,
then by Fact \ref{T334}(b) we have $ \ro(\rp ^n (\Col (\o _1 ,\k ^+)))\cong  \ro(\rp ^n (\ro(\CP _\k)))\cong  \ro(\rp ^n (\CP _\k))$
and (\ref{EQ324}) is true.
\hfill $\Box$
\paragraph{Boolean algebras invariant under reduced powers}
We will say that a Boolean algebra $\B$ is {\it invariant under reduced powers (IRP)}
iff $\ro (\rp ^n (\B ))\cong \ro (\B)$, for all $n\in \N$.
Using Fact \ref{T334}  it is easy to prove by induction that $\B$ is IRP iff $\ro (\rp (\B ))\cong \ro (\B)$.
$\B$ will be called {\it asymptotically IRP (AIRP)}
iff $\ro (\rp ^n (\B ))$ is an IRP algebra, for some $n\in \N$.
\begin{te}\label{T366}
If $(\ast)$ holds, then we have

(a) $\cf (\l)=\o _1$ $\Rightarrow$ $\Col (\l ,\k )$ is an IRP algebra $\Rightarrow$ $\cf (\l)\in [\o _1 ,\fh]$;

(b) $(\fh =\o _1)$ $\Col (\l ,\k )$ is an IRP algebra iff $\cf (\l)=\o _1$;

(c) $(\fh =\o _1 \lor (\fb =\fd \land 0^\sharp \mbox{ does not exist}))$
    If $\cf (\l)=\o$ then $\Col (\l ,\k )$ is an AIRP algebra;

(d) $(\fh =\o _1 \land  \mathrm{SCH})$
    If $\cf (\l)\leq\fc$ then $\Col (\l ,\k )$ is an AIRP algebra.
\end{te}
\dok
(a) If $\cf (\l)=\o _1$,
then by Fact \ref{T351} we have $\Col (\l ,\k )\cong \Col (\o _1 ,\k ^{<\l})$
and, by Theorem \ref{T340}(a), $\ro (\rp ^n (\Col (\o _1 ,\k ^{<\l}) ))\cong \Col (\o _1 ,\k ^{<\l})$.

We prove the contrapositive of the second implication.
If $\cf (\l)=\o$,
then by Fact \ref{T351} we have $\Col (\l ,\k )\cong \Col (\o ,\k ^{<\l})$
and, hence, the algebra $\Col (\l ,\k )$ collapses $\k ^{<\l}$ to $\o$,
while by Fact \ref{T300}(a) the algebra $\rp (\Col (\l ,\k ))$  preserves $\o _1$;
thus $\Col (\l ,\k )$ is not an IRP algebra.
If $\cf (\l)> \fh$, then by Fact \ref{T351} we have $\Col (\l ,\k )\cong \Col (\cf (\l) ,\k ^{<\l})$ and this algebra is $\fh$-distributive.
But by Theorem \ref{T352}(c) the algebra  $\rp (\Col (\l ,\k )$  is not $\fh$-distributive.

The rest follows from (a) and Theorems \ref{T335} and \ref{T364}.
\hfill $\Box$
\begin{te}\label{T369}
$(\fh =\o _1)$ A Boolean algebra $\B$ of size $\leq\fc$ is and IRP algebra iff $\ro (\B)\cong \Col (\o _1,\fc)$.
\end{te}
\dok
If $\ro (\rp (\B ))\cong \ro (\B)$,
then $\ro (\B )^+$ contains a dense subset isomorphic to $\rp (\B )^+$.
By Facts \ref{T300}(a) and \ref{T209} the poset $\rp (\B )^+$ is separative atomless $\s$-closed and collapses $\fc$ to $\o _1$.
Since $|\rp (\B )^+|\leq \fc ^\o=\fc$ we have $|\rp (\B )^+|=\fc$,
by Corollary \ref{T365} we have $\ro (\rp (\B )^+)\cong \Col (\o _1,\fc)$.
\hfill $\Box$
\paragraph{Acknowledgement.}
This research was supported by the Science Fund of the Republic of Serbia,
Program IDEAS, Grant No.\ 7750027:
{\it Set-theoretic, model-theoretic and Ramsey-theoretic
phenomena in mathematical structures: similarity and diversity}--SMART.

\end{document}